\documentclass[a4paper, 12pt,oneside,reqno]{amsart}
\usepackage[a4paper]{geometry}
\usepackage[matrix,arrow,curve,cmtip]{xy}

\usepackage{txfonts}
\DeclareMathAlphabet{\mathcal}{OMS}{cmsy}{m}{n} % but do not change mathcal symbols

\usepackage{
  amsmath,  % Contains many useful mathematical typesetting tools.
  amssymb,  % Contains many useful mathematical symbols.
  amsthm,   % Contains tools for theorem-like environments
  tikz,     % A versatile tool to draw diagrams
  fullpage, % More reasonable margins
  youngtab, % To draw young tableaux.
  ytableau, % this is the best young tableaux packege,
  thmtools, % Lets us define theorem-like environments easily.
  hyperref, % Better references.
  cite,     % Better citations.
  url       % URL handling
}

\usepackage{amscd}
\usepackage{euscript}
\usepackage{latexsym}
\usepackage[cp1251]{inputenc}
\usepackage[english]{babel}
\usepackage{mathrsfs}
\usepackage{graphicx}
\usepackage{keyval}
\usepackage{mathtools}
\usepackage{tikz}
\usetikzlibrary{arrows,shapes,snakes,automata,backgrounds,petri,through,positioning}
\usetikzlibrary{intersections}
\usepackage[symbol*]{footmisc}
\usepackage{verbatim}
\usepackage{tikz-cd}

\DeclareMathOperator{\B}{B}
\DeclareMathOperator{\C}{C}

\DeclareMathOperator{\dcobound}{d}

\DeclareMathOperator{\Homol}{H}

\DeclareMathOperator{\im}{Im}
\DeclareMathOperator{\id}{id}
\DeclareMathOperator{\Ker}{Ker}

\DeclareMathOperator{\Z}{Z}

\usepackage{etoolbox}
\patchcmd{\thebibliography}{\section*}{\paragraph}{}{}

\numberwithin{equation}{section}

\def\Vir{\mathrm{Vir}}
\def\Der{\mathop {\fam 0 Der} \nolimits}
\def\Cur{\mathop {\fam 0 Cur} \nolimits}
\def\End{\mathop {\fam 0 End} \nolimits}
\def\oo#1{\mathbin {{}_{(#1)}}}

\newtheorem{theorem}{Theorem}[section]
\newtheorem{lemma}[theorem]{Lemma}
\newtheorem{proposition}[theorem]{Proposition}
\newtheorem{corollary}[theorem]{Corollary}

\theoremstyle{definition}
\newtheorem{definition}[theorem]{Definition}
\newtheorem{example}[theorem]{Example}
\newtheorem{remark}[theorem]{Remark}

\title{Morse matching method for conformal cohomologies}
\author{H. Alhussein$^{1),2)}$, P.S. Kolesnikov$^{3),4)}$, 
V. Lopatkin$^{5)}$}

\address{$^{1)}$ Siberian State University of Telecommunication and Informatics, Novosibirsk, Russia}
\address{$^{2)}$Novosibirsk State University of Economics and Management, Novosibirsk, Russia}
\address{$^{3)}$Novosibirsk State University, Novosibirsk, Russia}
\address{$^{4)}$Sobolev Institute of Mathematics, Novosibirsk, Russia}
\address{$^{5)}$National Research University Higher School of Economics, Faculty of Computer Science,
Pokrovsky Boulevard 11, Moscow, 109028 Russia}

\subjclass[2020]{16E40; 17B68; 17B69; 16W25}
\keywords{Conformal algebra; Hochschild cohomology; Gr\"obner--Shirshov basis; Morse matching}

\def\Vir{\mathop {\fam 0 Vir} \nolimits}
\def\Cend{\mathop {\fam 0 Cend} \nolimits}
\def\oo#1{\mathbin {{}_{(#1)}} }

\begin{document}

\begin{abstract}
     We apply discrete algebraic Morse theory to the computation 
    of Hochschild cohomologies of associative conformal algebras. 
    As an example, we evaluate the dimensions 
    of the universal associative conformal envelope 
    $U(3)$ of the Virasoro Lie conformal algebra relative to 
    the associative locality $N=3$ on the generator with scalar coefficients.
\end{abstract}

\maketitle

\section{Introduction}

\subsection{Conformal algebras and their cohomologies.}\label{subsec1.1}

The notion of a conformal Lie algebra 
(also known as vertex Lie algebra) 
emerged in 
\cite{KacValgBeginners}
as an algebraic tool 
formalizing the properties of 
coefficients in the 
singular part of the operator 
product expansion (OPE) 
formula for chiral fields 
in 2-dimensional conformal field theory.
Namely, if $V$ is a vertex algebra,
$Y(a,z)$ and $Y(b,z)$ are  
vertex operators corresponding 
to $a,b\in V$
then the commutator $[Y(a,w),Y(b,z)]\in 
\End V[[z,z^{-1},w,w^{-1}]]$
may be expressed 
as a finite distribution with respect to 
derivatives of the formal delta-function:
\[
[Y(a,w),Y(b,z)] = \sum\limits_{s=0}^{N(a,b)-1}
 Y(c_s,z) \dfrac{1}{s!} \dfrac{\partial^s \delta(w-z)}{\partial z^s} ,
\]
$\delta(w-z) =
\sum\limits_{t\in \mathbb Z} w^tz^{-t-1}$.

The properties of binary algebraic operations 
$(\cdot \oo{s}\cdot)$
on 
$V$
given by
$(a \oo{s}b) = c_s$, 
$s\in \mathbb Z_+  =\{0,1,2,\dots \}$, $a,b\in V$,
along with the translation operator on $V$
lie in the background of the formal definition 
of a (Lie) conformal algebra. 

Lie conformal algebras also appear naturally (see \cite{Xu2000}) from 
the Hamiltonian formalism in the theory 
of differential equations of hydrodynamic type \cite{NovDubrovin}.
Namely, every Novikov algebra (or, more generally, a Gelfand--Dorfman 
algebra) gives rise to a conformal Lie algebra.
In particular, the simplest Poisson bracket on the phase space 
with one field function $v(x)$ given by 
\[
\{v(x),v(y)\} = 2v(x) \delta'(x-y) + v'(x) \delta(x-y)
\]
corresponds to the Virasoro conformal Lie algebra $\Vir $.

In a more general context, the family of operations 
$(\cdot \oo{s}\cdot)$, $s\in \mathbb Z_+$, 
may be defined 
on a space of pairwise local formal 
distributions over an arbitrary
(not necessarily Lie)
algebra. This leads to the definition 
of what is an associative (commutative, Jordan, etc.) 
conformal algebra, see \cite{Roitman99}.

A categorial approach to the theory of 
conformal algebras was proposed in \cite{BDK}. 
In this way, a conformal algebra is 
an algebra in the pseudo-tensor 
category of modules over the 
bialgebra $\mathbb C[\partial ]$ of polynomials in one variable
(with respect to the standard bialgebra structure).
This categorial definition also makes 
clear how to define cohomology of 
conformal algebras (c.~f. with \cite{BKV}).

The structure theory of Lie conformal algebras 
that are finite modules over $\mathbb C[\partial ]$
was established in \cite{DK1998}. In particular, 
every finite-dimensional simple Lie algebra $\mathfrak g$
gives rise to a simple conformal algebra $\Cur \mathfrak g$
embedded into the space of formal distributions 
over $\mathfrak g[t,t^{-1}]$,
and there is one exceptional 
simple Lie conformal algebra $\Vir $ 
mentioned above.

Irreducible modules over $\Cur \mathfrak g$ and $\Vir $
were described in \cite{ChengKac}, and the corresponding 
cohomologies were computed in \cite{BKV}.
Conformal cohomologies have the same relations 
to derivations, extensions, and deformations of 
conformal algebras as the ``ordinary'' ones. 

\subsection{Statement of the problem and main results}\label{subsec}
It is well-known (see, e.g., \cite[Ch.~XIII]{CartanEilenberg}) 
that for a Lie algebra $\mathfrak g$ acting on 
a $\mathfrak g$-module $V$ the cohomology groups 
$\Homol ^n(\mathfrak g, V)$
coincide with the Hochschild cohomology groups 
$\Homol ^n(U(\mathfrak g), V)$ 
of the universal associative enveloping 
algebra $U(\mathfrak g)$ with coefficients in the same~$V$ 
equipped with the induced structure of 
a $U(\mathfrak g)$-module.
The situation is different in the case for conformal algebras.

Given a conformal Lie algebra $L$ with a conformal $L$-module $M$, 
one may construct a series of universal enveloping 
associative conformal algebras corresponding to 
different associative locality functions on the 
generators \cite{Roit2000}. 
For example, consider the Virasoro conformal algebra $\Vir $.
One may fix a natural number $N$ 
and construct 
the associative conformal algebra $U(N)$ 
generated by a single 
element $v$ such that $(v\oo{n} v) = 0$ for $n\ge N$, and 
the commutation relations of $\Vir $ hold.
Obviously, $U(1)=0$; the algebra $U(2)$ 
is known as the Weyl conformal algebra 
(also denoted $\Cend_{1,x}$ \cite{BKL}). 
The structure of $U(3)=U(4)$ was studied in 
\cite{Kol2020-IAJC} by means of the Gr\"obner--Shirshov bases 
method.

It was shown in \cite{Kozlov2017} that the second Hochschild 
cohomology groups $\Homol ^2(U(2), M)$ are trivial for every conformal 
(bi-)module~$M$.

If $M$ is a conformal $\Vir$-module then it is not true in general 
that $M$ is a $U(N)$-module. A representation of $\Vir $ on $M$ 
is determined by the image $\rho(v)$ of $v$ in
the space of conformal endomorphisms $\Cend M$ (see \cite[Section~2.10]{KacValgBeginners}). 
For $M$ to be a $U(N)$-module, we need $\rho(v)\oo{n} \rho(v) = 0$
in $\Cend M$ for $n\ge N$.
The trivial 1-dimensional $\Vir $-module $M=\mathbb C$ 
is always a module over $U(N)$ since we have 
 $\rho(v)\oo{n} \rho(v) = 0$ for all $n\ge 0$. 
It is not hard to note (see \cite{AlKol-JMP}) 
that $\Homol ^n(U(2),\mathbb C) = 0$ for all $n\ge 1$.

Hence, even in the case of scalar coefficients there is no 
coincidence between $\Homol ^n(\Vir, \mathbb C)$ and 
$\Homol ^n(U(2), \mathbb C)$. The purpose of this note 
is to study the Hochschild cohomologies 
of $U(3)$, the next envelope in the series.
In \cite{AlKol-JMP}, it was found that 
$\dim \Homol^2(U(3),\mathbb C)=1$, but for 
higher Hochschild cohomologies the direct computation 
becomes too complicated since $U(3)$
is of quadratic growth (Gelfand--Kirillov dimension $=2$).

The purpose of this work is to develop 
a modification of the Morse matching method
for calculation of Hochschild cohomologies
of associative conformal algebras.
As an application, we find higher Hochschild cohomologies
of $U(3)$ with coefficients in~$\mathbb C$.

One more reason to study $U(3)$
rather than the smallest nonzero envelope $U(2)$
is the following.
If $M$ is a finite irreducible 
$\Vir $-module (one of those described in \cite{ChengKac}) 
then $\rho(v)\oo{n} \rho(v) = 0$ in $\Cend M$ for $n\ge 3$.
Hence, $U(3)$ is a more adequate associative conformal 
envelope of $\Vir $ than the Weyl conformal algebra~$U(2)$.

The calculations of conformal cohomologies in \cite{BDK}
is performed in an indirect way; they rest upon deep and nontrivial auxiliary construction. There is a natural question: whether one can arrive at these results in a more universal and natural fashion?

The aim of this paper is to develop a modification of the algebraic discrete Morse theory machinery (we call it the Morse matching method) to calculate cohomology of (associative) conformal algebras. We believe that calculation of homology of conformal algebras should be obtained in a natural manner; they should be deduced from an intrinsic structure of a combinatorial presentation of algebras (i.e., presentation via generators and relations). In this case this machinery looks natural and powerful and we demonstrate it in some examples.

\section{The basics of the algebraic discrete Morse theory}

In this section, we recall the basic definitions of 
the algebraic discrete Morse theory and show how to apply 
this theory to the computation of Anick resolutions.

\subsection{The Morse matching method}\label{subsec1.2}
Algebraic discrete Morse theory is an algebraic version of discrete Morse theory developed independently by Sk\"oldberg \cite{Sc}, 
and by J\"ollenbeck and Welker \cite{JW}. 
It allows us to construct, starting from a chain complex, a new
homotopically equivalent smaller complex using directed graphs. 
Here, for the
convenience of the reader, 
we present a short version of this machinery.
We follow
closely \cite[Chapter 2]{JW}, with minor simplifications and variations in 
notation.

Let $(\B_\bullet,\dcobound_\bullet)$ be a chain complex of free $R$-modules over a ring $R$,
\[
\B_0 \overset{\dcobound_1}{\leftarrow} \B_1 \overset{\dcobound_2}{\leftarrow} 
\B_2 
\overset{\dcobound_3}{\leftarrow } \dots .
\]

Suppose $X_n$, $n\ge 0$,  are some bases of the free $R$-modules $\B_n$. 
Consider the coordinates of the differentials 
$\dcobound_n: \B_n \to \B_{n-1}$ with respect to these bases:
\[
\dcobound_n(\mathbf c) = \sum\limits_{\mathbf c^\prime \in X_{n-1}}
 [\mathbf c: \mathbf c^\prime]\cdot \mathbf c^\prime ,
\]
where $\mathbf c \in X_n$, and $[\mathbf c:\mathbf c^\prime]$ are coefficients from $R$.

From these data, we construct a directed weighted graph 
$\Gamma(\B_\bullet) = (V,E)$. 
The
set of vertices $V$ of $\Gamma(\B_\bullet)$ 
is the union of all bases, $V = \bigcup_{n \ge 0} X_n$, 
and the set $E$ of weighted edges  $\mathbf c\to \mathbf c'$
consists of triples 
\[
\{ (\mathbf c, \mathbf c^\prime, [\mathbf c: \mathbf c^\prime])
\mid 
\mathbf c \in X_n, \, \mathbf c^\prime \in X_{n-1},\, [\mathbf c: \mathbf c^\prime] \ne 0\}.
\]
Here $\mathbf c$ and $\mathbf c^\prime $ are starting and ending vertices 
of an edge, $[\mathbf c:\mathbf c^\prime]\in R$ is its weight.

A subset $M \subseteq E$ of the set of edges is called an 
\emph{acyclic matching}, if it satisfies the following three conditions:

\begin{itemize}
    \item[(1)] (Matching) Each vertex $\mathbf v \in V$ lies in at most one edge in~$M$.
    
    \item[(2)] (Invertibility) For all edges 
    $(\mathbf c,\mathbf c', [\mathbf c:\mathbf c']) \in M$ 
     the weight $[\mathbf c: \mathbf c']$ lies in the center of 
     $R$ and is a unit (invertible) in $R$.
     
\item[(3)] (Acyclicity) The graph $\Gamma^M(\B_\bullet) = (V,E^M)$ 
constructed from the graph $\Gamma(\B_\bullet)$ has
no directed cycles, where
\[
E^M = (E \setminus M) \cup 
\{ (\mathbf c',\mathbf c, -[\mathbf c:\mathbf c']^{-1} \mid 
 (\mathbf c, \mathbf c', [\mathbf c:\mathbf c']) \in M\}.
\]
\end{itemize}

For an acyclic matching $M$ on the graph $\Gamma(\B_\bullet)$, we introduce the 
following notation:
\begin{itemize}
\item 
Define
\[
X^M_n = \{\mathbf c \in X_n \mid \mathbf c \text{ does not lie in any edge in } M \} .
\]
The vertices in $X^M_n$ are called 
{\em critical cells} of homological degree~$n$.
Let 
$\B^M_n$ be the subspace of $\B_n$ spanned by $X^M_n$.

%\item 
% Write $\mathbf c^\prime \le \mathbf c$  if $\mathbf c \in X_n$, $\mathbf c^\prime \in X_{n+1}$,
% and $[\mathbf c : \mathbf c^\prime] \ne 0$.

\item
$\mathscr{P}(\mathbf c,\mathbf c')$ is the set of paths from $\mathbf c$ to $\mathbf c'$ 
in $\Gamma^M(\B_\bullet)$.

\item 
The weight $\omega(\mathbf{p})$ of a path 
$\mathbf{p} = \mathbf c_1 \to \ldots \to \mathbf c_r \in \mathscr{P}(\mathbf c_1, \mathbf c_r)$
is defined as
the product of all weights of the edges $\mathbf c_k\to \mathbf c_{k+1}$, $k=1,\dots, r-1$.
Note that all these weights except one belong to the center of $R$, 
so the order of multiplication is not essential.
%\[
%\omega(\mathbf p) = \prod\limits_{k=1}^{r-1} \omega(c_k \to c_{k+1}),
%\quad
%\text{where}
%\quad
%\omega(c \to c^\prime) = \begin{cases} 
%- [c:c^\prime]^{-1} & \text{ if } c \le c^\prime \\ 
%\hskip 9pt [c:c^\prime]   & \text{ if } c^\prime \le c .
%\end{cases}
%\]

\end{itemize}

\begin{theorem}[{\cite[Theorem 2.2]{JW}}]%\label{Morse}
The chain complex $(\B_\bullet, \dcobound_\bullet)$ 
is homotopy equivalent to the complex 
$(\B^M_\bullet, \dcobound^M_\bullet)$, where 
the differential $\dcobound_n^M: \B^M_n \to \B^M_{n-1}$ is defined
as
\[
\dcobound_n^M(\mathbf c) = 
\sum\limits_{\mathbf c^\prime \in X_{n+1}^M} 
\Gamma_{\B^M_\bullet} (\mathbf c, \mathbf c'), 
\quad 
\Gamma_{\B^M_\bullet} (\mathbf c, \mathbf c') = 
\sum\limits_{\mathbf{p} \in \mathscr{P}(\mathbf c, \mathbf c^\prime)}
 \omega(\mathbf p) \mathbf c^\prime ,
\]
where $\mathbf c \in \B^M_n$.
\end{theorem}

\subsection{The Anick resolution for associative algebras.}
Let $\Lambda$ be an associative unital algebra over a field $\Bbbk$, 
the cokernel of the embedding map 
$\eta:\Bbbk \to \Lambda$ is denoted 
 $\Lambda/\Bbbk$. 
Assume further that $X$ 
is a set of generators of~$\Lambda$ (as of 
an associative algebra with identity). 
Suppose that a Gr\"obner--Shirshov basis 
of $\Lambda$ with respect to some 
ordering of the free monoid
$W(X)$  generated by $X$ is known.

As usual we denote by $\Lambda^e \coloneqq \Lambda 
\otimes \Lambda^{\mathrm{op}}$
the enveloping algebra for the algebra 
$\Lambda$. 
This algebra plays the role of $R$ in the previous subsection.
Following~\cite{JW} and \cite{Sc} we will see how to construct a free 
$\Lambda^e$-resolution for $\Lambda$.

Let us start with the two-sided bar resolution 
$\mathsf{B}_\bullet(\Lambda, \Lambda)$ 
which is a $\Lambda^e$-free resolution of $\Lambda$, 
where
\[
\mathsf{B}_n(\Lambda, \Lambda)\coloneqq \Lambda 
\otimes (\Lambda/\Bbbk)^{\otimes n} \otimes  \Lambda 
\cong \Lambda^e \otimes (\Lambda/\Bbbk)^{\otimes n}.
\]
The basis of $\mathsf B_n(\Lambda , \Lambda )$ 
as of $\Lambda ^e$-module 
is presented by 
$[a_1|a_2|\dots |a_n]$, where $a_i$ are nontrivial (i.e., not equal to 1) 
reduced normal forms relative to the fixed Gr\"obner--Shirshov basis.

The differential
$\mathrm d_n : \mathsf{B}_n(\Lambda, \Lambda)
   \to \mathsf{B}_{n-1}(\Lambda, \Lambda)$
is defined as follows:
\begin{multline*}
\mathrm{d}_n([a_1| \ldots | a_n]) = (a_1 \otimes 1) 
[a_2| \ldots |a_n]
+  \sum\limits_{i=1}^{n-1}(-1)^{i}   [a_1| \ldots | 
N(a_ia_{i+1})| \ldots | a_n]  \\
+ (-1)^n(1 \otimes a_n)[a_1| \ldots| a_{n-1}] .
\end{multline*}
Here $N(a_ia_{i+1})$ is the reduced normal form of the product $a_ia_{i+1}$.

In \cite{Anick1983}, Anick showed how to construct a free
resolution $\mathsf A_\bullet (\Lambda, \Lambda )$
which is essentially smaller than the bar-resolution but
homotopically equivalent to the latter. 
The linear bases $\Lambda ^{(n-1)}$ of 
$\mathsf A_n (\Lambda, \Lambda )$
consist of so called Anick chains related 
with the chosen Gr\"obner--Shirshov basis of~$\Lambda $.
The computation of a differential in the Anick resolution 
may be simplified by means of the Morse matching method.

For $w \in W(X)$, let $\Lambda_{w,p}$ be set of all the vertices 
$[w_1| \ldots | w_n]$ in 
$\Gamma(\mathsf{B}_\bullet(\Lambda,\Lambda))$ such that $w = w_1 \cdots w_n$ 
and~$p$ is the largest integer 
$p \ge -1$ for which $w_1\cdots w_{p+1} \in 
\Lambda^{(p)}$ is an Anick $p$-chain.
Let 
$\Lambda_w \coloneqq  \bigcup\limits_{p \ge -1}\Lambda_{w,p}$.

Define a partial matching $\mathcal{M}_w$ on 
$\Gamma(\mathsf{B}_\bullet(\Lambda,\Lambda))|_{\Lambda_w}$ by letting 
$\mathcal{M}_w$ consist of all edges
\[
[w_1| \ldots| w'_{p+2}|w''_{p+2}| \ldots| w_n] \to [w_1| \ldots| w_{p+2}| 
\ldots| w_m]
\]
where 
$w'_{p+2}w''_{p+2} = w_{p+2}$, 
$[w_1|  \ldots| w_m] \in \Lambda_{w,p}$, 
and $[w_1| \ldots| w_{p+1}|w'_{p+2}] \in \Lambda^{(p+1)}$ is an Anick 
$(p+1)$-chain.

\begin{theorem}[J\"ollenbeck--Sc\"oldberg--Welker]\label{JSW}
The set of edges $\mathcal{M} = \bigcup_{w}\mathcal{M}_w$ is a Morse matching on 
$\Gamma(\mathsf{B}_\bullet(\Lambda,\Lambda))$, 
the critical cells of 
homological degree $p\ge 0$ are exactly $(p-1)$th Anick chains 
$\Lambda^{(p-1)}$.
\end{theorem}

For the matching $\mathcal{M} = 
\bigcup\limits_{\omega}\mathcal{M}_\omega$ we have the following result.

\begin{proposition}[{\cite[Chapter 5]{JW}, \cite[Lemma 9 and Theorem 5]{Sc}}]\label{HAR}
The set of edges $\mathcal{M} = \bigcup\limits_{\omega}\mathcal{M}_\omega$ is a 
Morse matching on $\Gamma(\mathsf{B}_\bullet(\Lambda, \Lambda))$, with Anick 
chains as critical cells. Moreover, the complex 
$\mathsf{A}_\bullet(\Lambda,\Lambda)= \mathsf B_\bullet ^{\mathcal M} (\Lambda, \Lambda )$
which is defined as 
follows:
\[
\mathsf{A}_{n+1}(\Lambda,\Lambda) = \Lambda^e \otimes \Bbbk\Lambda^{(n)}, 
\qquad
\mathrm{d}_{n+1}(v) = \sum\limits_{\mathbf v' \in \Lambda^{(n-1)}}
\Gamma_{\mathsf B^{\mathcal M}_\bullet }(\mathbf v,\mathbf v')\mathbf v',
\]
is a free $\Lambda^e$ resolution of $\Lambda$.
\end{proposition}

The Morse matching method provides us 
a powerful tool 
to compute differentials in the Anick resolution 
$\mathsf A_\bullet(\Lambda , \Lambda)$ and in the complex 
$\mathsf A_\bullet 
= \mathsf A_\bullet(\Lambda , \Lambda) \otimes_{\Lambda^e} \Bbbk $.

\begin{example}
 Let $\Lambda = U(\mathfrak g)$ be 
the universal enveloping associative algebra of a Lie algebra $\mathfrak g$
with an ordered basis $X$. Then $\{ xy-yx-[x,y] \mid  x,y\in X,\, x>y\}$
is a Gr\"obner--Shirshov basis of $\Lambda $ 
with respect to the deg-lex ordering,
and the corresponding Anick resolution is 
exactly the Chevalley--Eilenberg complex of $\mathfrak g$. 
\end{example}

\begin{example}\label{exmp:U(2)-coeff}
Let $\Lambda $
be the unital associative algebra 
over a field $\Bbbk$ generated by
infinite family of elements $v(n)$, $n\ge 0$,
relative to the defining relations 
$v(n)v(m)-v(m)v(n)-(n-m)v(n+m-1)$ and
$v(n+2)v(m)-2v(n+1)v(m+1)+v(n)v(m+2)$, 
$n,m\ge 0$.
The Gr\"obner--Shirshov basis consists of 
relations $v(n)v(m) -v(0)v(n+m)-nv(n+m-1)$, $n\ge 1$, $m\ge 0$.
\end{example}

Then $\Lambda^{(n)} = \{[v(p_1)|\dots |v(p_{n})|v(p_{n+1})] \mid 
p_1,\dots, p_n\ge 1 \}$.

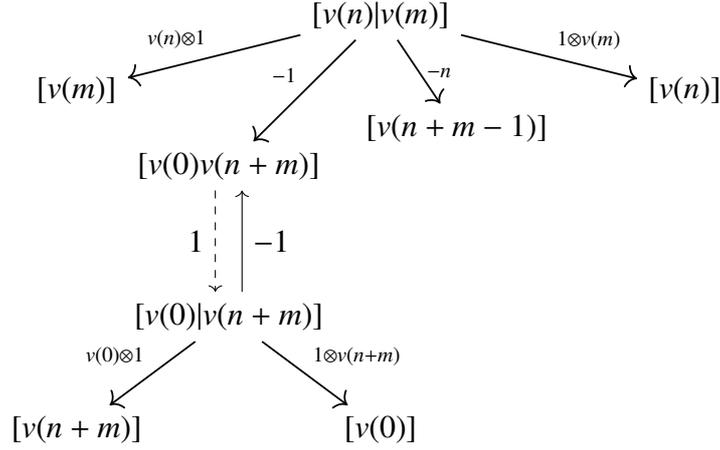
\begin{figure} 
\centering
\begin{tikzpicture}[commutative diagrams/every diagram]
\node (A1) at (0,0) {$[v(n)|v(m)]$};

\node (B1) at (-4,-1) {$[v(m)]$};
\node (B2) at (-2, -2) {$[v(0)v(n+m)]$};
\node (B3) at (1,-1.5)  {$[v(n+m-1)]$};
\node (B4) at (4,-1)  {$[v(n)]$};

\node (C1) at (-2,-4) {$[v(0)|v(n+m)]$} ;

\node (D1) at (-4, -5.5) {$[v(n+m)]$};
\node (D2) at (0,-5.5) {$[v(0)]$};

\path[commutative diagrams/.cd,every arrow, every label]
  (A1) edge node [swap] {$v(n) \otimes 1$} (B1)
  (A1) edge node [swap] {$-1$} (B2)
  (A1) edge node [right] {$-n$} (B3)
  (A1) edge node [above right]{$1 \otimes v(m)$} (B4)
  (C1) edge node [swap] {$v(0)\otimes 1$} (D1)
  (C1) edge node [above right] {$1\otimes v(n+m)$} (D2)
   ;

 \path[->]
  ([xshift=5pt]C1.north) edge node[right] {$-1$} ([xshift=5pt]B2.south);
  \path[->,dashed]
  ([xshift=-5pt]B2.south) edge node[left]{$1$} ([xshift=-5pt]C1.north);

\end{tikzpicture}
    \caption{Calculation of $\dcobound_2$ in Example~\ref{exmp:U(2)-coeff}}
    \label{fig:Fig0}
\end{figure}

In order to find $\dcobound_2(\mathbf v)$ for $\mathbf v = [v(n)|v(m)])$
we have to apply the formulas from Proposition~\ref{HAR}
with $\Gamma(\mathbf v,\mathbf v')$ calculated via the graph 
on Figure~\ref{fig:Fig0}:
\[
\dcobound_2([v(n)|v(m)])
=v(n)\otimes [v(m)] - v(0)\otimes [v(n+m)] - [v(0)]\otimes v(n+m)
-n[v(n+m-1)] + [v(n)]\otimes v(m).
\]

The calculation of $\dcobound_3([v(n)|v(m)|v(p)])
\in\mathsf A_2(\Lambda, \Lambda)$
by means of Proposition~\ref{HAR} leads to the following formula:
\begin{multline}\label{eq:Diff-U(2)}
\dcobound_3([v(n)|v(m)|v(p)]) = v(n)\otimes [v(m)|v(p)] - v(0)\otimes [v(n+m)|v(p)]
-n[v(n+m-1)|v(p)] \\
+ m[v(n)|v(m+p-1)] +[v(n)|v(0)]\otimes v(m+p) 
- [v(n)|v(m)]\otimes v(p) \\
+ n [v(n-1)|v(m+p)] +v(0)\otimes [v(n)|v(m+p)].
\end{multline}
Indeed, the corresponding segment of 
$\Gamma ^{\mathcal M}(\mathsf B_\bullet (\Lambda,\Lambda))$
is stated on Figure~\ref{fig:Fig1}.
The edges from $\mathcal M$ in the initial graph 
are drawn dashed, the vertices corresponding to the Anick chains are boxed.
They are critical cells corresponding to $\mathcal M$, 
but the vertices of the form $[v(0)|v(k)]$ are not critical cells 
since they lie on the matching edges of the ``smaller level'', see Figure~\ref{fig:Fig2}. This is why such vertices are ``dead ends'' of 
the calculation process, they do not contribute to the expression 
for~$\dcobound_n$.

In order to find the value of $\dcobound_n$ on a particular 
Anick chain we do not need to expand vertices
in $(\Lambda/\Bbbk)^{\otimes n}$ that contain 
two or more components $v(0)$ since neither 
of them may produce an Anick chain from 
$(\Lambda/\Bbbk)^{\otimes (n-1)}$. 
The same observation will be used in Section~\ref{Sec4}.

In the following, we will omit the tensor sign to make
expressions shorter. 

\begin{figure}
    \centering
     \begin{tikzpicture}[commutative diagrams/every diagram]
  \node (m p)       at (-4,0.5)  {$\boxed{[v(m)|v(p)]}$ };
  \node (n m)       at (4,0.5)   {$ \boxed{[v(n)|v(m)]} $};
  \node (n m p)     at (0,0)     {$[v(n)|v(m)| v(p)]$ };
  \node (on+m p)    at (-5,-0.5) {$[v(0) v(n+m)|v(p)]$};
  \node (n om+p)    at (5,-0.5)  {$[v(n)|v(0)v(m+p)]$};
  \node (n+m-1 p)   at (-1.7,-2) {$\boxed{[v(n+m-1)|v(p)]}$};
  \node (n m+p-1)   at (1.7,-2)  {$\boxed{[v(n)|v(m+p-1)]}$};
  \node (o n+m p)   at (-5,-3)   {$[v(0)|v(n+m)|v(p)]$};
  \node (n+m pl)    at (-8,-4)   {$\boxed{[v(n+m)|v(p)]}$};
  \node (o n+m+p-1) at (-7, -5)  {$[v(0)|v(n+m+p-1)]$};
  \node (o on+m+p)  at (-4, -6)  {$[v(0)|v(0)v(n+m+p)]$};
  \node (n+m pr)    at (-3, -4)  {$[v(0)|v(n+m)]$};
  \node (o o n+m+p) at (-4,-8)   {$[v(0)|v(0)|v(n+m+p)]$};
  \node (o n+m+p)   at (-7,-9)   {$[v(0)|v(n+m+p)]$};
  \node (oo n+m+p)  at (-4,-10)  {$[v(0)v(0)|v(n+m+p)]$};
  \node (o o)       at (-1, -9)  {$[v(0)|v(0)]$};
  \node (n o m+p)   at (5,-3)    {$[v(n)|v(0)|v(m+p)]$};
  \node (o m+p)     at (1, -4)   {$[v(0)|v(m+p)]$};
  \node (on m+p)    at (2,-6)    {$[v(0)v(n)|v(m+p)]$};
  \node (n-1 m+p)   at (5, -5)   {$\boxed{[v(n-1)|v(m+p)]}$};
  \node (n o)       at (6,-4)    {$\boxed{[v(n)|v(0)]}$};
  \node (o n m+p)   at (2, -8)   {$[v(0)|v(n)|v(m+p)]$};
  \node (o n)       at (6, -7)   {$[v(0)|v(n)]$};
  \node (o n+m+p-1l)at (6, -9)   {$[v(0)|v(n+m+p-1)]$};
  \node (o on+m+pl) at (5, -10)  {$[v(0)|v(0)v(n+m+p)]$};
  \node (n m+p)     at (0.5,-10) {$\boxed{[v(n)|v(m+p)]}$}; 
  \node (o o n+m+pl)at (5, -12)  {$[v(0)|v(0)|v(n+m+p)]$};
  \node (o n+m+pl)   at (0, -13) {$[v(0)|v(n+m+p)]$};
  \node (oo n+m+pl)  at (3, -14){$[v(0)v(0)|v(n+m+p)]$};
  \node (o ol)       at(6, -13.5) {$[v(0)|v(0)]$};
  
  \path[commutative diagrams/.cd,every arrow, every label]
  (o o n+m+pl) edge node[swap] {$v(0)\otimes 1$} (o n+m+pl)
  (o o n+m+pl) edge node[swap] {$-1$} (oo n+m+pl)
  (o o n+m+pl) edge node[right] {$-1 \otimes v(n+m+p)$} (o ol)
  ;

  \path[->, dashed]
  ([xshift=5pt]o o n+m+pl.north) edge node[right] {$1$} ([xshift=5pt]o on+m+pl.south);
  \path[->]
  ([xshift=-5pt] o on+m+pl.south) edge node[left]{$-1$} ([xshift=-5pt]o o n+m+pl.north); 
  
  \path[commutative diagrams/.cd,every arrow, every label]
  (o n m+p) edge node[swap] {$v(0)\otimes 1$} (n m+p)
  (o n m+p) edge node[swap] {$1$} (o on+m+pl)
  (o n m+p) edge node[above] {$n$} (o n+m+p-1l)
  (o n m+p) edge node [swap] {$-1 \otimes v(m+p)$} (o n)
  ;

  \path[->, dashed]
  ([xshift=5pt]o n m+p.north) edge node[right] {$-1$} ([xshift=5pt]on m+p.south);
  \path[->]
  ([xshift=-5pt] on m+p.south) edge node[left]{$1$} ([xshift=-5pt]o n m+p.north);

  \path[commutative diagrams/.cd,every arrow, every label]
  (n o m+p) edge node[swap] {$v(n) \otimes 1$} (o m+p)
  (n o m+p) edge node[right] {$-1 \otimes v(m+p)$} (n o)
  (n o m+p) edge node[swap] {$-1$} (on m+p)
  (n o m+p) edge node[swap] {$-n$} (n-1 m+p)
  ;

  \path[commutative diagrams/.cd,every arrow, every label]
  (n m p) edge node[above] {$v(n) \otimes 1$} (m p)
  (n m p) edge node[above] {$-1 \otimes v(p)$} (n m)
  (n m p) edge node[above] {$-1$} (on+m p)
  (n m p) edge node[above] {$1$} (n om+p)
  (n m p) edge node[swap] {$-n$} (n+m-1 p)
  (n m p) edge node[above right] {$m$} (n m+p-1)
   ;
  
 %arrows from [v(0)v(n+m)|v(p)]--->[v(0)|v(n+m)|v(p)]
  \path[->, dashed]
  ([xshift=5pt]o n+m p.north) edge node[right] {$-1$} ([xshift=5pt]on+m p.south);
  \path[->]
  ([xshift=-5pt]on+m p.south) edge node[left]{$1$} ([xshift=-5pt]o n+m p.north);

 \path[->, dashed]
  ([xshift=5pt]o o n+m+p.north) edge node[right] {$1$} ([xshift=5pt]o on+m+p.south);
  \path[->]
  ([xshift=-5pt]o on+m+p.south) edge node[left]{$-1$} ([xshift=-5pt]o o n+m+p.north);

 \path[commutative diagrams/.cd,every arrow, every label]
  (o n+m p) edge node [yshift=-2pt, xshift = -3pt, swap] {$v(0) \otimes 1$} (n+m pl)
  (o n+m p) edge node [swap] {$1$} (o n+m+p-1)
  (o n+m p) edge node [yshift=-12pt, xshift = 2pt] {$n+m$} (o on+m+p)
  (o n+m p) edge node [yshift=-4pt] {$-1\otimes v(p)$} (n+m pr)
  ;
  
 \path[commutative diagrams/.cd,every arrow, every label]
  (o o n+m+p) edge node [swap] {$v(0) \otimes 1$} (o n+m+p)
  (o o n+m+p) edge node [swap] {$-1$} (oo n+m+p)
  (o o n+m+p) edge node [xshift = 4pt, yshift = -1pt] {$-1\otimes v(n+m+p)$} (o o)
   ;
  
  \path[->, dashed]
  ([xshift=5pt]n o m+p.north) edge node[right] {$1$} ([xshift=5pt]n om+p.south);
  \path[->]
  ([xshift=-5pt]n om+p.south) edge node[left]{$-1$} ([xshift=-5pt]n o m+p.north);
  
\end{tikzpicture}
    \caption{Calculation of $\dcobound_3$ in Example~\ref{exmp:U(2)-coeff}}
    \label{fig:Fig1}
\end{figure}

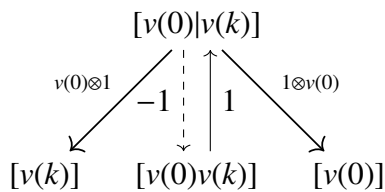
\begin{figure}
    \centering
\begin{tikzpicture}[commutative diagrams/every diagram]
\node (A) at (0,0) {$[v(0)|v(k)]$};
\node (B) at (-2,-2) {$[v(k)]$};
\node (C) at (0, -2) {$[v(0)v(k)]$};
\node (D) at (2,-2)  {$[v(0)]$};

\path[commutative diagrams/.cd,every arrow, every label]
  (A) edge node [swap] {$v(0) \otimes 1$} (B)
  (A) edge node [above right] {$1 \otimes v(0)$} (D)
   ;

 \path[->]
  ([xshift=5pt]C.north) edge node[right] {$1$} ([xshift=5pt]A.south);
  \path[->,dashed]
  ([xshift=-5pt]A.south) edge node[left]{$-1$} ([xshift=-5pt]C.north);

\end{tikzpicture}
    \caption{A ``dead end'' segment of the Morse matching graph in Example~\ref{exmp:U(2)-coeff}}
    \label{fig:Fig2}
\end{figure}

%-----------------------------

In this paper, we will adjust the Morse mathcing method 
for associative conformal algebras and, as an application,
calculate all Hochschild cohomology 
groups for $U(3)$ with coefficients in the scalar 
module $M=\Bbbk $. 
Note that $\Homol ^2(U(3),\Bbbk )$
was earlier found 
in \cite{AlKol-JMP} by ``elementary'' methods
which do not work even for the third cohomology 
since the computations become too bulky. 

%With the Morse matching method, we show that 
%\[
%\dim \Homol ^2(U(3), \Bbbk ) = \dim \Homol ^3(U(3), \Bbbk ) = 1
%\]
%and $\Homol ^n(U(3),\Bbbk ) = 0$ for $n\ge 4$ (and for $n=1$),
%which is in good correspondence with 
%$\Homol ^n(\Vir , \Bbbk )$.

\section{Conformal algebras}

Let $\Bbbk $ be a field of characteristic zero, 
and let $\mathbb Z_+$ stand for the set of nonnegative 
integers.

\subsection{Definition and examples of conformal algebras}\label{2.1} 

\begin{definition}[\cite{KacValgBeginners}]
A conformal algebra is a linear space $C$ equipped with 
a linear operator $\partial : C\to C$
and with a polynomial-valued map ($\lambda $-product)
\[
(\cdot \oo\lambda \cdot ): C\otimes C \to 
\Bbbk [\partial,\lambda]\otimes_{\Bbbk[\partial ]} C\cong 
\Bbbk[\lambda ]\otimes C,
\]
where $\lambda $ is a formal variable, 
satisfying the following axioms:
\begin{gather}
    (\partial a\oo\lambda b) = -\lambda (a\oo\lambda b), 
      \label{eq:3/2-lin(1)}\\
    (a\oo\lambda \partial b) = (\partial+\lambda) (a\oo\lambda b),
      \label{eq:3/2-lin(2)}
\end{gather}
for all $a,b\in C$.
\end{definition}

The coefficients of the polynomial $(a\oo\lambda b)$ at 
$\lambda^n/n!$ are  denoted $a\oo{n} b$, $n\in \mathbb Z_+$.
The definition of a conformal algebra may be stated 
equivalently in terms of multiple binary operations 
$(\cdot \oo{n} \cdot)$, for all non-negative 
integers~$n$.

Note that the torsion of $C$ as of $\Bbbk [\partial]$-module has the annihilating property:
$\mathrm{tor}\,C \oo\lambda C = C \oo\lambda \mathrm{tor}\,C = 0$.

For every conformal algebra $C$ one may construct 
an orinary algebra $A=\mathcal A(C)$ in the following way. 
As a linear space, 
\[
A = \Bbbk [t,t^{-1}] \otimes _{\Bbbk [\partial ]} C,
\]
where $\partial $ acts on $\Bbbk[t,t^{-1}]$ as $- d/dt$.
Denote $t^n\otimes_{\Bbbk [\partial ]} a$, $a\in C$, $n\in \mathbb Z$, 
by $a(n)$.
The operation on $A$ is given by 
\[
a(n)\cdot b(m) = \sum\limits_{s\ge 0} \binom{n}{s} (a\oo{s} b) (n+m-s).
\]
The map 
\begin{equation}\label{eq:Coeff-Derivation}
\partial : a(n)\mapsto n a(n-1)
\end{equation}
is a well-defined derivation of $A$.

It is easy to see that the linear span of all $a(n)$ with 
$n\in \mathbb Z_+$ is a $\partial$-invariant subalgebra 
of $A$ denoted $\mathcal A_+(C)$. This algebra 
plays an important role in the cohomology theory 
of conformal algebras.

The algebra $A=\mathcal A(C)$ constructed has the following properties 
\cite{Roitman99}. 
First, $C$ embeds into a 
formal distribution conformal algebra over $A$, i.e., the map $\imath : C\to A[[z,z^{-1}]]$, 
$a\mapsto \sum\limits_{n\in \mathbb Z} a(n) z^{-n-1}$, 
is injective, $\imath (\partial a) = d\imath (a)/dz$,
and
\[
\imath (a\oo\lambda b)(z) = \mathop{\fam0 Res}\limits_{w=0}
\big( 
\exp^{\lambda (w-z)} \imath(a)(w)\cdot \imath (b)(z)\big ).
\]
Second, the algebra $A$ is universal among all ordinary 
algebras $B$ such that $C$ maps to a formal distribution 
conformal algebra over~$B$.

The correspondence $C\mapsto \mathcal A(C)$ is functorial and 
provides a foundation for a definition of what is a 
Lie or associative conformal algebra. 
A conformal algebra $C$ is said to be associative (Lie, etc.) 
if so is $\mathcal A(C)$. 
In terms of $\lambda $-products, the associativity may 
be expressed as 
\begin{equation}\label{eq:Conf_assoc}
    (a\oo{\lambda } (b\oo\mu c)) = ((a\oo{\lambda } b)\oo{\mu+\lambda } c)\in C[\lambda, \mu ], 
    \quad a,b,c\in C.
\end{equation}
The anti-commutativity and Jacobi identity are equivalent to 
\begin{gather}
     \label{eq:Conf-A-comm}
    (a\oo\lambda b) = - (b\oo{-\partial-\lambda } a), \\
    \label{eq:Conf-Jacobi}
(a\oo{\lambda } (b\oo\mu c)) - (b\oo\mu (a\oo{\lambda } c)) 
= ((a\oo{\lambda } b)\oo{\mu+\lambda } c), 
\end{gather}
respectively.

An associative 
conformal algebra $C$ turns into a Lie conformal algebra $C^{(-)}$
under a new $\lambda $-product corresponding 
to the ordinary commutator in $\mathcal A(C)$:
\[
[a\oo\lambda b ] = (a\oo\lambda b) - (b\oo{-\partial -\lambda } a),
\quad a,b\in C.
\]

\begin{remark}
An associative (resp., Lie) conformal algebra 
may be considered as an associative (resp., Lie) 
algebra in a certain pseudo-tensor category related 
with $\Bbbk [\partial ]$, as proposed in \cite{BDK}. 
This categorial approach to the definition 
of a variety of conformal algebras 
is equivalent to the described above.
\end{remark}

The functor $(-)$ from the category of associative 
conformal algebras to the category of Lie conformal algebras 
has no left adjoint one. 
Indeed, given a Lie conformal algebra $L$ with a
$\lambda $-product 
$[\cdot \oo\lambda \cdot]$, one may construct 
a series of associative conformal algebras $C$ such that 
$L\subset C^{(-)}$, but there is no upper bound for the function 
$N_C:L\times L\to \mathbb Z_+$,
$N_C(a,b) = \deg_\lambda (a\oo\lambda b)+1$.
So the universal  enveloping associative conformal algebra 
for a Lie conformal algebra does not exist in general.
However, if we set up additional conditions to restrict 
the function $N_C$ then we may construct a ``partial'' 
universal envelope. 

\begin{example}
Let $A$ be an ordinary algebra over $\Bbbk $. 
Then the free $\Bbbk [\partial]$-module 
$C = \Bbbk [\partial] \otimes A$ equipped with 
\[
(\partial^n \otimes a) \oo\lambda (\partial ^m \otimes b)
 = -\lambda ^n (\lambda +\partial )^m \otimes ab, 
 \quad a,b\in A,
\]
is a conformal algebra denoted $\Cur A$, the current conformal 
algebra over $A$.
\end{example}

Clearly, it is enough to define the $\lambda $-product 
only on the generators of a free $\Bbbk [\partial]$-module 
due to \eqref{eq:3/2-lin(1)}, \eqref{eq:3/2-lin(2)}.

If $A$ is associative or Lie then so is $\Cur A$ as a conformal 
algebra: $\mathcal A_+ (\Cur A) = A[t,t^{-1}]$, 
$\mathcal A_+(\Cur A) = A[t]$.

\begin{example}
Let $A$ be an associative algebra with a locally 
nilpotent derivation $D: A\to A$. 
Then the free $\Bbbk [\partial]$-module 
$C = \Bbbk [\partial] \otimes A$ equipped with 
\[
(1\otimes a) \oo\lambda  (1\otimes b) = 
  \sum\limits_{s\ge 0} \dfrac{\lambda^s}{s!}(1 \otimes aD^s(b)), 
 \quad a,b\in A,
\]
is an associative conformal algebra denoted 
$\mathop{\fam 0 Diff}(A,D)$.
\end{example}

% \? Check:
% \mathcal A(Diff (A,D)) is a skew differential polynomial algebra
% with id automorphism and derivation D

If $A=M_n(\Bbbk [x])$ with $D=d/dx$ then 
$\mathop{\fam 0 Diff}(A,D)$
is denoted $\Cend_n$, the algebra of conformal endomorphisms
of the free $\Bbbk [\partial]$-module of rank~$n$.

If $A=x\Bbbk[x]$ is the augmentation ideal of the polynomial 
algebra and $D = d/dx$ as above then 
$\mathop{\fam 0 Diff}(A,D)$
is known as the Weyl conformal algebra $\Cend_{1,x}$.

\begin{example}
The rank 1 free $\Bbbk [\partial]$-module 
$V = \Bbbk [\partial ]v$
is a Lie conformal algebra with respect to the $\lambda$-product 
\[
[v\oo\lambda v] = (\partial +2\lambda )v.
\]
This structure is known as the Virasoro conformal algebra $\Vir $, the exceptional simple finite Lie conformal algebra
\cite{DK1998}.
\end{example}

The coefficient algebra for $\Vir $ is the Witt algebra of 
vector fields on a circle: 
$\mathcal A(\Vir) = \Der \Bbbk [t,t^{-1}]$, 
$\mathcal A_+(\Vir) = \Der \Bbbk [t]$.

Let $\mathfrak g$ be a Lie algebra.
If $U(\mathfrak g)_0$ stands for the augmentation 
ideal of the universal enveloping associative algebra 
$U(\mathfrak g)$ then $\Cur U(\mathfrak g)_0$
is the universal enveloping associative 
conformal algebra of $L=\Cur \mathfrak g$ 
in the class of all those associative 
conformal envelopes $C$ of $L$ for which 
$N_C(\mathfrak g, \mathfrak g)\le 1$. Let 
us denote this conformal algebra by $U(\Cur\mathfrak g; N=1)$.
The structures of $U(\Cur\mathfrak g; N=2)$
and $U(\Cur\mathfrak g; N=3)$ are more complicated, 
they were studied in \cite{KolKoz_AlgRep} 
by means of the Gr\"obner--Shirshov bases technique
which was developed 
for associative conformal algebras in~\cite{BFK}.

Given an integer $N\ge 2$, 
one may construct an associative conformal envelope
$U(N)$
for the Virasoro Lie conformal algebra $\Vir $ with 
a generator $v$ 
which is universal
in the class of all such envelopes $C$ that $N_C(v,v)\le N$.
For example, $U(2)= U(\Vir ; N=2)$ is the Weyl conformal 
algebra~$\Cend_{1,x}$; the structure of $U(3) = U(\Vir ; N=3)$ 
is more complicated, it was studied in 
\cite{Kol_ProcICAC}. 
For the associative conformal algebra $U(2)$, the algebra 
$\mathcal A_+(U(2))$ is considered in Example~\ref{exmp:U(2)-coeff}.

The definition of a (bi-)module over an associative (or Lie) 
conformal algebra is very natural \cite{ChengKac}. 
In brief, a bimodule $M$ over an associative 
conformal algebra $C$ is a linear space equipped with 
a linear operator $\partial $ (denoted by the same symbol
as in $C$) and with left and right $\lambda $-actions
\[
C\otimes M \to M[\lambda ],\quad
M\otimes C \to M[\lambda ],
\]
satisfying the analogues of \eqref{eq:3/2-lin(1)}, \eqref{eq:3/2-lin(2)}, and \eqref{eq:Conf_assoc}. 
Obviously, these conditions are equivalent to the statement 
that $C\oplus M$ equipped with 
\[
(a+u)\oo\lambda (b+v) = (a\oo\lambda b) + (a\oo\lambda v + u\oo\lambda b), \quad a,b\in C,\ u,v\in M,
\]
is an associative conformal algebra (split null extension 
of $C$ by means of $M$).
In a similar way, conformal modules over 
Lie conformal algebras are defined.

If $M$ is a conformal bimodule over an associative 
conformal algebra $C$ then the same space $M$ may 
be considered as a bimodule over the 
ordinary associative algebra $\mathcal A_+(C)$. 
Namely, 
\[
a(n)u = a\oo{n} u, \quad ua(n) = \{u\oo{n} a\}, \quad u\in M, \ a\in C,\ n\ge 0,
\]
where 
\[
\{u\oo{n} a\} = (-1)^n\sum\limits_{s\ge 0 } \frac{(-\partial)^s}{s!} (u\oo{n+s} a).
\]

\begin{example}\label{exmp:Vir-modules}
Let $L=\Vir$, $M = \Bbbk [\partial ]u$ 
is a 1-generated free $\Bbbk [\partial]$-module 
equipped with a left action of $\Vir $ via
\[
(v\oo\lambda u) = (\alpha +\partial +\Delta\lambda )u, 
\]
where $\alpha, \Delta \in \Bbbk $. 
These are all irreducible (for $\Delta \ne 0$) 
conformal modules over $\Vir $ (see \cite{ChengKac})
denoted $M(\alpha, \Delta)$.
\end{example}

For every finite conformal $\Vir$-module there 
exists a composition series of submodules 
in which all quotients are either $M(\alpha, \Delta)$
or just 1-dimensional $\Bbbk $ considered as a trivial 
module \cite{ChengKac}. 

The action of $\Vir $ on $M(\alpha, \Delta )$
corresponds to the following homomorphism of conformal algebras 
(representation)
\[
\begin{aligned}
\rho :{} & \Vir \to (\Cend_1)^{(-)} = \mathrm{gc}_1, \\
& v\mapsto f=x-\Delta \partial +\alpha .
\end{aligned}
\]
It is easy to see that if $\Delta \ne 0$ then 
$\deg (f\oo\lambda f)=2$, so for $C=\Cend_1$ we have 
$N_C(\rho(v), \rho(v)) = 3$. 
Therefore, the structure of a conformal 
$\Vir$-module on $M(\alpha, \Delta)$ 
may not be extended to $U(2)$ but may be extended to $U(3)$. 
This is a motivation to study $U(3)$ 
rather than a more simple Weyl conformal algebra. 
One more reason is related with the cohomology theory 
of these algebras. 

\subsection{Conformal cohomologies}\label{2.2} 

The study of cohomologies for conformal algebras was initiated in
\cite{BKV}. 
Let us state the main definitions and results 
concerning the Hochschild cohomologies of associative conformal 
algebras. We will focus on the reduced complex \cite{BDK} 
since its cohomologies have the expected relations 
to derivations, extensions, and deformations of 
associative conformal algebras. Moreover, 
the reduced complex coincides with the 
construction arising when one considers 
conformal algebras and their modules 
in the framework of pseudo-tensor categories \cite{BDK}, 
see also \cite{KolKoz_CMP}.

Let $C$ be an associative conformal algebra, 
and let $M$ be a conformal bimodule over~$C$.
The Hochschild complex $\C^\bullet (C,M)$
consists of the cochain spaces $\C^n (C,M)$, $n=1,2,\ldots$,
each of them is the space of all maps 
\[
\varphi_{\bar \lambda }: 
C^{\otimes n}\to M[\lambda_1,\ldots, \lambda_{n-1}],
\]
where $\bar\lambda = (\lambda_1,\ldots, \lambda_{n-1})$,
satisfying the analogues of \eqref{eq:3/2-lin(1)} and 
\eqref{eq:3/2-lin(2)}: 
\[
\begin{gathered}
 \varphi_{\bar \lambda }(a_1,\ldots,\partial a_i, \ldots , a_n) =
 -\lambda_i \varphi _{\bar \lambda }(a_1,\ldots, a_n),\quad i=1,\ldots, n-1, \\
 \varphi_{\bar \lambda }(a_1,\ldots, \partial a_n) = 
 (\partial +\lambda_1 + \cdots + \lambda_{n-1})\varphi_{\bar \lambda }(a_1,\ldots, a_n).
\end{gathered}
\]
The Hochschild differential 
\[
\dcobound_n : \C^n(C,M) \to \C^{n+1}(C,M)
\]
is given by
\begin{multline}\nonumber
(\dcobound_n\varphi)_{\bar\lambda }
(a_1,\ldots, a_{n+1}) = 
a_1 \oo{\lambda_1} 
\varphi_{\bar\lambda_0} (a_2, \ldots , a_{n+1}) 
 + \sum\limits_{i = 1}^{n} (-1)^i\varphi_{\bar\lambda_i}
  (a_1,\ldots , a_i \oo{\lambda_i} a_{i+1},\ldots , a_{n+1}) \\
  +  (-1)^{n+1}\varphi_{\bar \lambda_{n+1} }(a_1 ,\dots , a_n) \oo{\lambda_1+\cdots +\lambda_n} a_{n+1},
\end{multline}
where 
$\bar \lambda = (\lambda_1, \ldots, \lambda_{n})$, 
$\bar\lambda_0 = (\lambda_2, \ldots, \lambda_{n})$, 
$\bar\lambda_i = (\lambda_1, \ldots, \lambda_i+\lambda_{i+1},
\ldots , \lambda_n )$, 
$\bar \lambda_{n+1} = (\lambda_1, \ldots, \lambda_{n-1})$.

It is common to complete the complex $\C^\bullet (C,M)$ described above with 
$\C^0(C,M) = M/\partial M$ and $\dcobound_0: \C^0(C,M)\to \C^1(C,M)$, where 
\[
\dcobound_0(u+\partial M): a\mapsto \{a\oo{0} u\} - u\oo{0} a.
\]

The elements of the second cohomology group $\Homol^2(C,M) = 
\Ker \dcobound_2/\im \dcobound_1 $ are in one-to-one 
correspondence with 
the classes of equivalent null extensions $E$, 
\[
0\to M\to E\to C\to 0, \quad (M\oo\lambda M)=0.
\]
The similar statement holds for the cohomologies 
of Lie conformal algebras \cite{BKV} and their central 
extensions. Let us state two examples to demonstrate 
this relation.

\begin{example}
The Virasoro Lie conformal algebra $\Vir $ 
has a unique (up to a scalar multiple) central extension 
\[
0\to \Bbbk \to \Vir_c \to \Vir \to 0,\quad c\in \Bbbk , 
\]
where 
\[
[v\oo\lambda v] = (2\lambda+\partial) v+ \lambda^3 c \in \Vir_c = \Vir\oplus \Bbbk .
\]
So $\Homol ^2(\Vir ,\Bbbk )$ is 1-dimensional.
\end{example}

\begin{example}
Let $\mathfrak g$ be a finite-dimensional semisimple Lie algebra 
with the Killing form $\langle \cdot,\cdot \rangle $. 
Then $\Cur\mathfrak g$ has a nontrivial central
extension
$\Cur\mathfrak g \oplus \Bbbk $, where 
\[
[a\oo\lambda b] = ab + \lambda \langle a,b\rangle , 
\quad a,b\in \mathfrak g.
\]
\end{example}

The Hochschild complex $\C^\bullet  (C,M)$
for an associative conformal algebra $C$ and a conformal
$C$-bimodule $M$ may also be constructed as follows. 
Let us consider $M$ as a bimodule over the
associative algebra $\Lambda =\mathcal A_+(C)\oplus \Bbbk 1$
(with external identity).
Both $M$ and $ \mathcal A_+(C) $
carry linear operators (both denoted $\partial $) 
such that their sum is a derivation of the split null 
extension $ \mathcal A_+(C) \oplus M$. Then define
\[
\partial_n^* : \C^n( \mathcal A_+(C),M) \to \C^n( \mathcal A_+(C) ,M)
\]
as follows:
\[
(\partial_n^* f)(\alpha_1,\ldots, \alpha_n)
= \partial f(\alpha_1,\ldots, \alpha_n)
+ \sum\limits_{i=1}^n 
f(\alpha_1,\ldots, \partial\alpha_i, \ldots, \alpha_n),
\]
where $\partial (a(n)) = n a(n-1)$ for $a\in C$, $n\ge 0$.
The maps $\partial^*_n $ form 
a morphism of complexes
$\partial^*_\bullet: \C^\bullet (\mathcal A_+(C),M)
\to 
\C^\bullet (\mathcal A_+(C),M)$, and the 
following statement holds.

\begin{proposition}[{\cite[Theorem 6.1]{BKV}}]\label{prop:MainTool}
$\C^\bullet (C,M)\cong 
\C^\bullet ( \mathcal A_+(C), M)/\partial^*_\bullet \C^\bullet ( \mathcal A_+(C),M)$.
\end{proposition}

The latter complex is called a reduced Hochschild complex of an associative conformal algebra $C$
with coefficients in~$M$.

\subsection{Reduced complex and Hochschild cohomologies 
of current conformal algebras}\label{subsec:ReducedComplex}

Let $\Lambda = \mathcal A_+(C)\oplus \Bbbk 1$,
with $\varepsilon (a(n))=0$ for all $a\in C$, $n\ge 0$. 
Assume $\Lambda $ acts trivially on the 1-dimensional space 
$\Bbbk $, i.e., $\alpha 1\beta = \varepsilon(\alpha \beta )$  for 
$\alpha,\beta \in \Lambda $.
In order to calculate Hochschild cohomologies 
of $\Lambda $ it is enough to apply the Hom functor 
to the complex
$\mathsf{B}_\bullet 
= \mathsf{B}_\bullet (\Lambda, \Lambda)\otimes_{\Lambda^e}\Bbbk $, 
 where 
$\mathsf{B}_\bullet (\Lambda, \Lambda)$ is the 
two-sided bar resolution for~$\Lambda $, 
$\Lambda ^e$ stands for the 
enveloping algebra $\Lambda \otimes \Lambda^{op}$. 
Apparently, 
\[
\mathsf B_{n} = (\Lambda /\Bbbk )^{\otimes n}
\]
Denote by $\dcobound_n: \mathsf{B}_{n} \to \mathsf{B}_{n-1}$
the differential of $\mathsf{B}_\bullet$
induced by the bar differential.
The dual map $\dcobound_n^*$ is the Hochschild 
differential $\C^{n-1} \to \C^{n}$, 
$\C^i = \C^i(\mathcal A_+(C),\Bbbk )$.

Recall that $\Lambda $ is equipped with a derivation $\partial $
such that $\partial (a(n)) = na(n-1)$ for $a\in C$, $n\ge 0$, 
and the 1-dimensional bimodule carries trivial derivation.
Let us extend $\partial $ to the $\Lambda^e$-linear map 
$\partial_n : \mathsf B_{n} \to \mathsf B_{n}$
by the rule 
\begin{equation}\label{eq:Partial_on_B}
\partial_n : [\lambda_1 | \ldots | \lambda_n] \mapsto 
 \sum\limits_{i=1}^n [\lambda_1 | \ldots | \partial(\lambda_i)| \ldots | \lambda_n].
\end{equation}
Then the dual map $\partial_n^*$ is the morphism mentioned 
in Proposition~\ref{prop:MainTool}. 
Since $\C^n/\partial_n^*\C^n \simeq (\Ker \partial_n)^*$
by the Fredholm principle, 
we reduce the problem of computing conformal cohomologies
$\Homol^n(C,\Bbbk )$, $n\ge 1$, 
to the application of the Hom functor to the complex
\begin{equation}\label{eq:KerComplex}
\cdots \leftarrow \Ker\partial_n \leftarrow \Ker\partial_{n+1}\leftarrow \cdots . 
\end{equation}
%
% Ker \partial \embedding_to  B_n
% dual: \C^n \onto (\Ker\partial)^*
% kernel = orth.complement to Ker\partial = image of \partial^*
The arrows here are restrictions of $\dcobound_i$
onto $\Ker\partial_i \subset \mathsf B_{i}$.
%Finding $\dim \Homol^n(C,\Bbbk )$ is now a 
%linear algebra problem.

\begin{example}
Let $A$ be an associative algebra, 
and let $C=\Cur A$ be the current conformal algebra over $\Bbbk $. 
Then $\Lambda /\Bbbk  = \mathcal A_+(C) = A[t]$
is spanned by $a(m)$, $a\in A$, $m\ge 0$, 
so that 
\[
a(n)b(m) = (ab)(n+m).
\]
\end{example}

Note that for the current conformal algebra 
the kernel of $\partial _n$ is easy to find. 
There is an isomorphism of linear spaces
\[
\mathsf B_n= A[t]^{\otimes n}
\simeq \Bbbk [x_1,\ldots , x_n]
\otimes A^{\otimes n}, 
\]
where 
$[a_1(m_1)| \dots |a_n(m_n)]$ corresponds to 
$x_1^{m_1}\cdots x_n^{m_n}\otimes (a_1|\dots |a_n)$. 
The map $\partial _n$ acts on the space  
$\Bbbk [x_1,\ldots , x_n]
\otimes A^{\otimes n}$
as the differential operator
\[
D =\dfrac{\partial}{\partial x_1} +\cdots + \dfrac{\partial}{\partial x_n} 
\]
Thus we may identify $\Ker \partial_n $
with the subspace 
$K_n\otimes A^{\otimes n}$
where $K_n$ is the kernel of $D$ in $\Bbbk [x_1,\dots, x_n]$.
It is not hard to note that $K_n$ is the subalgebra generated 
by $y_i = x_{i+1}-x_i$, $i=1,\ldots, n-1$.

With these notations, a value of $\dcobound_n$ on $\Ker\partial_n$
may be expressed as follows:
\[
\dcobound_n (f\otimes \mathbf v) = 
\sum\limits_{i=1}^{n-1} \dcobound_n^{(i)}(f)
\otimes \mathbf v^{(i)} ,\quad f\in K_n, \ \mathbf v\in A^{\otimes n},
\]
where 
\[
\dcobound_n^{(i)} (f(y_1,\ldots, y_{n-1})
= (-1)^i f(y_1,\ldots, y_{i-1}, 0, y_i,\ldots, y_{n-2}), 
\]
and 
$\mathbf v^{(i)} = a_1\otimes \cdots \otimes a_ia_{i+1}\otimes \cdots \otimes a_n$
for 
$\mathbf v=a_1\otimes \cdots \otimes a_n$.

In general, we may deduce the following statement on the Hochschild cohomologies 
of current associative conformal algebras 
with trivial coefficients.

\begin{theorem}\label{thm:CurrentCohomologies}
Let $A$ be an associative algebra acting 
trivially on $\Bbbk $, and let 
$\Lambda  = A\oplus \Bbbk $
be the augmented algebra with $\varepsilon(A)=0$.
Then 
$\Homol ^1(\Cur A, \Bbbk ) 
= \Homol^1(\Lambda ,\Bbbk ) 
= (A/A^2)^*$.
For every $n\ge 2$, if 
$\dim \Homol^1(\Lambda  ,\Bbbk ) = \cdots 
=\dim \Homol^{n-1}(\Lambda  ,\Bbbk ) = 0$
then 
\[
\Homol^{n}(\Cur A,\Bbbk ) = 
\Homol^n(\Lambda ,\Bbbk ).
\]
\end{theorem}

\begin{proof}
Let $\delta_n :  A^{\otimes n}\to A^{\otimes (n-1)} $
stand for the derived differential of the bar resolution of $\Lambda $, 
i.e., $\Homol^n(\Lambda ,\Bbbk ) = \Ker \delta_{n+1}^*/\im \delta_{n}^*$.

Let us start with $n=1$.
By definition, $\Ker\partial_1$ 
consists of $[x(0)]=1\otimes x\in K_n\otimes A$, 
$x\in A$, 
and it is easy to see that $\Ker \partial_2$
is spanned by 
\[
e_m(a,b) = y_1^m\otimes (a\otimes b)=
\sum\limits_{s=0}^m (-1)^s\binom{m}{s} 
[a(m-s)|b(s)], \quad 
a,b\in A,\ m\ge 0.
\]
Obviously, 
\[
\dcobound_2 e_m(a,b) = \begin{cases}
[ab(0)], & m=0, \\
0, & m>0.
\end{cases}
\]
Therefore,
$\Ker \dcobound_1 / \im \dcobound_2$
in the complex \eqref{eq:KerComplex}
is isomorphic to $A/A^2$.
Hence, 
$\Homol ^1(C, \Bbbk ) 
= \Homol^1(\Lambda ,\Bbbk ) 
= (A/A^2)^*$. 

Proceed to $n=2$. Note that 
\[
\dcobound_2\left ( 
 \sum\limits_{n\ge 0} y_1^n\otimes \mathbf b_n
\right) = -\mathbf b_0^{(1)}, 
\quad 
\dcobound_3\left ( 
 \sum\limits_{n,m\ge 0} y_1^n y_2^m\otimes \mathbf c_{n,m}
\right) =
\sum\limits_{n\ge 0} y_1^n \big (
\mathbf c_{n,0}^{(2)} - \mathbf c_{0,n}^{(1)}\big).
\]
Suppose $\mathbf u \in \Ker \dcobound_2$, 
$\mathbf u = \sum\limits_{n\ge 0} y_1^n\otimes \mathbf b_n$.
Since $A=A^2$, for every $n>0$ we may find $\mathbf c_n \in A^{\otimes 3}$ 
such that $\mathbf c_n^{(1)}= \mathbf b_n$.
Then for 
$\mathbf v=\sum\limits_{n\ge 0} y_2^n\otimes \mathbf c_n $
we have 
$\mathbf u-\dcobound_3(\mathbf v) \in 1\otimes A^{\otimes 2}$.
On the other hand, 
$\Ker\dcobound_2 \cap (1\otimes A^{\otimes 2})
=1\otimes \Ker\delta_2$,
as well as 
$\im\dcobound_3 \cap (1\otimes A^{\otimes 2})
=1\otimes \im\delta_3$.
Therefore, 
\[
\Ker \dcobound_2 /\im \dcobound_3 \simeq 
\Ker \delta_2 /\im \delta_3 , 
\]
as desired.

Consider the general case $n\ge 2$. Suppose 
$\mathbf u\in \Ker \dcobound_n$,
$\mathbf u\notin 1\otimes A^{\otimes n}$. 
Let us present $\mathbf u$ as
\[
\mathbf u  = y_1^{m_1}\cdots y_{n-1}^{m_{n-1}}\otimes \mathbf b_0 +\dots, 
\]
where the first (``leading'') term is chosen in
the following way. 
Find maximal $i$ such that a positive power of $y_i$ appears in $\mathbf u$.
Among all such monomials (with $y_i$), choose the maximal one in the lexicographic sense.

If $i=n-1$ then find $\mathbf c\in A^{\otimes (n+1)}$
such that $\mathbf c^{(n)} = \mathbf b_0$, and note that 
for 
$\mathbf v = y_1^{m_1}\cdots y_{n-1}^{m_{n-1}}\otimes \mathbf c$
we have
\[
\dcobound_{n+1}(\mathbf v)= 
(-1)^n
y_1^{m_1}\cdots y_{n-1}^{m_{n-1}}\otimes \mathbf c^{(n)}
+ \dots ,
\]
where all remaining terms do not contain $y_{n-1}$.
Hence, $\mathbf u-(-1)^n\dcobound_{n+1}(\mathbf v)$
has a smaller 
leading term than~$\mathbf u$.

If $i<n-1$ then consider the coefficient at
$y_1^{m_1}\cdots y_i^{m_i}$ in $\dcobound_n(\mathbf u)$.
By the choice of the leading term in $\mathbf u$, we have 
\[
\dcobound_n(\mathbf u)
=y_1^{m_1}\cdots y_i^{m_i} \otimes 
\big ( 
(-1)^{i+1}\mathbf b_0^{(i+1)} +(-1)^{i+2} \mathbf b_0^{(i+2)} + \cdots 
+ (-1)^{n-1} \mathbf b_0^{(n-1)}\big) + \dots,  
\]
where all remaining terms either do not contain $y_i$
or lexicographically smaller than 
$y_1^{m_1}\cdots y_i^{m_i}$.
Therefore, 
\[
(\id^{\otimes i} \otimes \delta_{n-i})(\mathbf b_0)=0
\]
and since $\Homol^{n-i}(\Lambda ,\Bbbk )=0$ we may find 
$c\in A^{\otimes (n+1)}$ such that 
\[
(\id^{\otimes i} \otimes \delta_{n+1-i})(\mathbf c)= \mathbf b_0.
\]
As in the previous case, we may reduce the 
leading term of $u$ by considering 
$\mathbf u - \dcobound_{n+1}(\mathbf v)$
for 
\[
\mathbf v = (-1)^{i} y_1^{m_1}\cdots y_i^{m_i} \otimes \mathbf c.
\]

The reduction of the leading term described above 
shows that every class in 
$\Ker \dcobound_n/\im \dcobound_{n+1}$
contains an element from $1\otimes A^{\otimes n}$.
The rest of the proof is completely similar to $n=2$
case.
\end{proof}

\begin{corollary}
If $A=M_n(\Bbbk )$ then 
$H^n(\Cur A, \Bbbk ) = 0$
for all $n\ge 1$.
\end{corollary}

% M_n(A) has H^1 trivial in ANY bimodule
% Hochschild => H^n = 0
% now apply theorem
%

The result obtained for 
Hochschild cohomologies of 
current associative conformal algebras 
looks different from 
what was proved in \cite{BKV} for current 
Lie conformal algebras:
if $\mathfrak g$ is a semisimple finite-dimensional 
Lie conformal algebra 
then 
$\Homol^n(\Cur \mathfrak g, \Bbbk )=
\Homol^n(\mathfrak g, \Bbbk )+\Homol^{n+1} 
(\mathfrak g, \Bbbk )$ for all $n\ge 0$. 

In particular, if $A = \mathfrak gU(\mathfrak g)$
then $\Cur A$ is the universal enveloping 
associative conformal algebra 
for $\Cur \mathfrak g$ relative to the locality 
bound $N=1$ on the elements of~$\mathfrak g$, 
$\Lambda = U(\mathfrak g)$.
If $[\mathfrak g,\mathfrak g] = \mathfrak g$ then 
$A^2=A$, so by Theorem \ref{thm:CurrentCohomologies}
$\Homol ^2(\Cur A,\Bbbk ) =
\Homol^2(\Lambda , \Bbbk ) 
= \Homol^2(\mathfrak g, \Bbbk ) 
= 0$ in contrast to the 
1-dimensional 
Lie conformal cohomology group
$\Homol^2(\Cur \mathfrak g, \Bbbk )$.

A similar picture appears when we consider 
cohomology groups with trivial coefficients 
for the Virasoro Lie conformal algebra $\Vir $
and its universal enveloping associative 
conformal algebras $U(N)$, $N=2,3,\dots $.
It was proved in \cite{Kozlov2017} that 
$\Homol ^2(U(2),M)=0$ for every conformal bimodule 
over the Weyl conformal algebra $U(2)$ in contrast 
to 1-dimensional $\Homol ^2(\Vir, \Bbbk )$. 
This is the reason to study $\Homol^n(U(3),\Bbbk )$
which is the aim of Section~\ref{Sec4}.

\subsection{On the Anick resolution for differential algebras}

In the previous section we exploited the 
bar resolution 
$\mathsf{B}_\bullet (\Lambda, \Lambda)$ 
and the complex $\mathsf{B}_\bullet \otimes_{\Lambda^e}\Bbbk $
for the augmented 
associative algebra 
$\Lambda = \mathcal A_+(C)\oplus \Bbbk $
to compute Hochschild cohomology groups 
of a conformal algebra $C$ in the case 
when $C$ is the current associative 
conformal algebra. In more complicated cases, 
the computation with bar resolution 
becomes much harder, so it is reasonable 
to replace the bar resolution with a more 
compact Anick (two-sided) resolution 
$\mathsf A_\bullet (\Lambda,\Lambda )$. 
The first problem in this route 
is to translate the mapping $\partial_\bullet $ from 
\eqref{eq:Partial_on_B} to the complex 
$\mathsf A_\bullet = 
\mathsf A_\bullet (\Lambda,\Lambda )\otimes_{\Lambda^e}\Bbbk $. 

%Let $\Lambda $ be an augmented algebra 
%acting trivially on a bimodule $M=\Bbbk $. 
%The latter means, as above, that 
%$\lambda u \mu = \varepsilon(\lambda \mu)u$
%for $\lambda, \mu \in \Lambda $, $u\in M$.
%Assume  $\partial $ is a derivation of 
%$\Lambda $ such that $\varepsilon(\partial(\Lambda ))=0$.

The mapping 
$\partial_n: \mathsf{B}_n\to \mathsf{B}_n $
defined by \eqref{eq:Partial_on_B}
is a morphism of complexes, 
$\dcobound_n\circ \partial_n = \partial_{n-1}\circ \dcobound_n$.
(This is not the case when we consider a nontrivial conformal module.)
Suppose 
$\mathsf{A}_\bullet = 
\mathsf{A}_\bullet(\Lambda,\Lambda)\otimes _{\Lambda^e} \Bbbk $,
$\mathsf{A}_n = \Lambda^{(n-1)}$. 

The space $\mathsf B_n$ is spanned by elements of the form 
$[u_1|\dots |u_n]$, where $u_i$ are nontrivial 
reduced words in the generators of $\Lambda $.
Then $\mathsf A_n$ is a subspace of $\mathsf B_n$
spanned by the Anick $(n-1)$-chains.
Let us define the linear projection 
$\pi_{\mathsf A_n}: \mathsf{B}_n \to \mathsf{A}_n$
assuming 
\[
 \pi_{\mathsf A_n}(\mathbf a)=
 \begin{cases}
  0, & \text{if $\mathbf a$ is not an Anick chain}, \\
  \mathbf a, & \text{otherwise}.
\end{cases}
\]

\begin{proposition}
 Let $\partial_n: \mathsf{B_n} \to \mathsf{B}_n$ 
 be defined by \eqref{eq:Partial_on_B} then in terms of the complex 
 $\mathsf{A}_\bullet$ we have
  \[
  \tilde \partial (\mathbf a_n) = \sum_{\mathbf b_n \in \mathsf{B}_n} 
  \Gamma_{\mathsf{A}_\bullet}(\mathbf a_n,\mathbf b_n) \cdot \pi_{\mathsf{A}_n}
  (\partial_n(\mathbf{b}_n)).
 \]
\end{proposition}

\begin{proof}
Since the Anick complex can be obtained from the bar complex by Morse matching machinery then by \cite[Appendix B, (B2), (B3)]{JW} for any $n>0$ we have the following commutative diagram
 \[
  \xymatrix{
   \mathsf{A}_n \ar@{->}[r]^{\tilde \partial_n} \ar@{->}[d]_{\mathrm{g}_n} & \mathsf{A}_n \\
   \mathsf{B}_n \ar@{->}[r]_{\partial_n} & \mathsf{B}_n \ar@{->}[u]_{\mathrm{f}_n}
  }
 \]
 where the maps $\mathrm{g}_\bullet$, $\mathrm{f}_\bullet$ are defined as follows
 \[
  \mathrm{f}_n(\mathbf{b}_n) = \sum_{\mathbf{a}_n\in \mathsf{A}_n}
     \Gamma_{\mathsf{B}_\bullet}(\mathbf{b}_n,\mathbf{a}_n) \mathbf{a}_n, 
  \quad
  \mathrm{g}_n(\mathbf{a}_n) = \sum_{\mathbf{b}_n \in \mathsf{B}_n}
     \Gamma_{\mathsf{A}_\bullet}(\mathbf{a}_n,\mathbf{b}_n)\mathbf{b}_n.
\]
By Lemma \cite[Appendix B, Lemma B.3]{JW} these maps define a chain homotopy between 
the resolutions $\mathsf{A}_\bullet$ and $\mathsf{B}_\bullet$. 
Therefore we can define $\tilde \partial$ as follows:
$\tilde \partial (\mathbf{a}_n) = 
 (\mathrm{f}_n\circ \partial_n \circ \mathrm{g}_n) (\mathbf{a}_n)$ 
 for every $n >0$.
We thus have
 \begin{equation}\label{eq:DerivTilda}
  \tilde \partial (\mathbf{a}_n) = \sum_{\mathbf{b}_n \in \mathsf{B}_n}
  \sum_{\mathbf{a}_n' \in \mathsf{A}_n} \Gamma_{\mathsf{A}_\bullet}(\mathbf{a}_n,\mathbf{b}_n)
   \Gamma_{\mathsf{B}_\bullet}(\partial_n(\mathbf{b}_n), \mathbf{a}'_n))\cdot \mathbf{a}_n'.
 \end{equation}
On the other hand,  
\[
 \Gamma_{\mathsf{B}_\bullet}(\mathbf{b}_n, \mathbf{a}_n)= \begin{cases}
  1 & \mbox{if $\mathbf{b}_n = \mathbf{a}_n \in \mathsf{A}_n$,} \\
  0 & \mbox{otherwise,}
 \end{cases}
\] 
and the statement follows.
\end{proof}

\begin{example}
Let $\mathfrak g$ be a Lie algebra with a 
linearly ordered basis $X$. Then consider 
\[
A = 
\Bbbk \langle a(n), a\in X, n\ge 0 \mid 
a(n)b(m)-b(m)a(n)-[a,b](n+m) \rangle .
\]
Then $\partial : a(n)\mapsto n a(n-1)$
is a derivation of $A$. The corresponding 
augmented algebra $\Lambda = A\oplus \Bbbk 1$
with $\varepsilon(a(n))=0$ is just 
the universal enveloping 
associative algebra of $\mathfrak g[t]$.
\end{example}

Let us order the generators of $A$ as follows:
\[
a(n)>b(m) \iff n>m \text{ or }n=m\text{ and } a>b.
\]
Then $[a(1)|b(0)]\in \Lambda ^{(1)}$ for $a<b$.
It is not hard to construct the graph $\Gamma({{\mathsf A}_\bullet})$
(similar to Example~\ref{exmp:U(2)-coeff}) 
and find
$\mathrm{g}_2([a(1)|b(0)]) = [a(1)|b(0)]-[b(0)|a(1)]$.
Since $[a(0)|b(0)]$ is not a chain for $a<b$, 
we obtain 
\[
\tilde\partial_2 ([a(1)|b(0)]) = - [b(0)|a(0)].
\]

\begin{example}
Let $\Lambda $ be the algebra from Example \ref{exmp:U(2)-coeff}
acting trivially on the scalar module $M=\Bbbk $. 
Note that $\partial : v(n)\mapsto n v(n-1)$
defines a derivation on $\Lambda $.
We will denote an element of the form 
$[v(n_1)|\ldots |v(n_k)]$ by 
$[n_1|\ldots |n_k]$ for brevity. 

Consider the complex $\mathsf A_\bullet 
= \mathsf A_\bullet(\Lambda ,\Lambda)\otimes_{\Lambda^e} \Bbbk $.
Then $\dcobound_3([2|1|0]) = 2[1|1] - [2|0]$
and 
$\dcobound_3([1|1|0]) = 0$
by \eqref{eq:Diff-U(2)}.
Following Figure~\ref{fig:Fig1}, we may see that 
$\mathrm g_3([2|1|0]) = [2|1|0]-[0|3|0]+3[0|0|3] - [2|0|1]
+ [0|2|1] - [0|0|3]$. Hence, 
\[
\tilde\partial_3([2|1|0]) = 2[1|1|0],
\]
and by \eqref{eq:Diff-U(2)}
we have $\dcobound_3([1|1|0])=0$.
In a similar way, one may calculate
$\tilde\partial_2(2[1|1] - [2|0]) = 2[1|0]-2[1|0] =0$
in compliance with $\dcobound_3\tilde \partial_3 = \tilde\partial_2\dcobound_3$.
\end{example}

Therefore, in order 
to calculate conformal cohomologies $H^n(C,\Bbbk )$
following the scheme of 
Section~\ref{subsec:ReducedComplex},
we have to study the complex 
\[
\dots \leftarrow \Ker \tilde\partial_{n-1} \leftarrow \Ker \tilde \partial_n \leftarrow \dots, 
\]
where the arrows are restrictions of the Anick 
differential $\dcobound_n: \mathsf A_n\to \mathsf A_{n-1}$
to the kernel of $\tilde\partial_n$.

\section{Hochschild cohomologies of the $N=3$ universal associative envelope
of the Virasoro conformal algebra}\label{Sec4}

%\subsection{}\label{4.1}
%Weyl conformal algebra $N=2$: check $H^n=0$

\subsection{Gr\"obner--Shirshov basis of $\mathcal A_+(U(3))$}\label{4.2}

By definition, the universal enveloping 
associative conformal algebra $U(3)$ 
of the Virasoro Lie conformal algebra $\Vir $
relative to the locality bound $N=3$ 
is generated by a single element $v$
such that $v\oo{n} v = 0$ for $n\ge 3$. 
The remaining defining relation of $U(3)$ is 
\begin{equation}\label{eq:vir-comm}
2v\oo{1} v - \partial (v\oo{2} v)= 2v.
\end{equation}
The algebra $\mathcal A_+(U(3))$ 
is generated by the elements $v(n)$, $n\ge 0$, 
relative to the following relations (see \cite{Roit2000}):
\begin{gather}
 v(n)v(m) -3v(n-1)v(m+1) + 3v(n-2)v(m+2) - v(n-3)v(m+3) = 0,\quad n\ge 3, \ m\ge 0, 
 \label{eq:u3-defn-loc}
 \\
v(n)v(m) - v(m)v(n) = (n-m)v(n+m-1), \quad n>m\ge 0. 
\label{eq:u3-defn-comm}
\end{gather}

Let us fix the deg-lex order on the set 
of words of the form $v(n_1)\dots v(n_k)$
assuming $v(n)>v(m)$ iff $n>m$.

\begin{theorem}\label{thm:GSB-U(3)}
The Gr\"obner--Shirshov basis of $\mathcal A_+(U(3))$
consists of the relations
\begin{equation}\label{eq:u3-gsb-10}
v(1)v(0)=v(0)v(1)+v(0),
\end{equation}
\begin{multline}\label{eq:u3-gsb}
    v(n)v(m) =\frac{nm}{n+m-1} v(1)v(n+m-1)-\frac{(n-1)(m-1)}{n+m-1}v(0)v(n+m) \\
    +\frac{n(n-1)}{n+m-1}v(n+m-1),\quad n\ge 2.
\end{multline}
\end{theorem}

\begin{proof}
First, let us prove that \eqref{eq:u3-gsb}
hold on $U(3)$ for all $n\ge 2$ and $m\ge 0$. For $m=0,1$ and for every $n\ge 2$ 
this is just \eqref{eq:u3-defn-comm}.
Proceed to the case $n=2$. 
Calculate the $(m+1)$th Fourier coefficient of 
\eqref{eq:vir-comm}:
\[
2(v(1)v(m+1) - v(0)v(m+2)) + (m+1) (v(2)v(m)
-2v(1)v(m+1) + v(0)v(m+2) = 2v(m+1),
\]
i.e., \eqref{eq:u3-gsb} holds for $n=2$
and $m\ge 0$. It remains to apply 
induction on $n\ge 2$ using 
\eqref{eq:u3-defn-comm} for the induction 
step. For example, in the generic case $n\ge 5$
we have 
\begin{multline*}
v(n)v(m) 
= 3v(n-1)v(m+1) -3v(n-2)v(m+2) + v(n-3)v(m+3) \\
=\dfrac{1}{n+m-1}
\big (
3(n-1)(m+1) v(1)v(n+m-1) - 3(n-2)m v(0)v(n+m) 
  + 3(n-1)(n-2) v(n+m-1) \\
-3(n-2)(m+2) v(1)v(n+m-1) + 3(n-3)(m+1) v(0)v(n+m) 
  - 3(n-2)(n-3) v(n+m-1) \\
+(n-3)(m+3) v(1)v(n+m-1) - (n-4)(m+2) v(0)v(n+m) 
  + (n-3)(n-4) v(n+m-1) 
\big ) \\
=\dfrac{1}{n+m-1}
\big(
   nm v(1)v(n+m-1)
 -(n-1)(m-1) v(0)v(n+m)
 + n(n-1) v(n+m-1)
\big ) ,
\end{multline*}
as desired. (For $n=3,4$ the only difference is 
that we should not expand some terms via \eqref{eq:u3-gsb}.)

Next, let us make sure that the set of relations
\eqref{eq:u3-gsb-10}, \eqref{eq:u3-gsb}
is closed under composition. 
This may be done in a straightforward way, but we may simply note that the reduced words are linearly 
independent in $\mathcal A_+(U(3))$. 
Indeed, the words reduced modulo \eqref{eq:u3-gsb-10}, \eqref{eq:u3-gsb}
are of the form 
\[
v_{k,p;m} = v(0)^k v(1)^pv(m), 
\]
where $n\ge 1$ for $p>0$ or $n\ge 0$ for $p=0$.
On the other hand, the basis of $U(3)$ as of an 
$H$-module was found in \cite{Kol_ProcICAC}, 
it consists of 
\[
v_k = (v\oo{0})^k v,
\quad
v_{k,p}=(v\oo0 )^k (v \oo1)^p (v\oo2 v)
\]
for $k,p\ge 0$. Here $v\oo{n}$ stands for the 
operator of $n$th conformal  multiplication 
$(v\oo{n}\cdot )$ on $U(3)$.
It is easy to calculate the principal terms 
of the Fourier coefficients for $v_k$, $v_{k,p}$:
\[
v_{k}(m) = v_{k,0;m}, 
\quad 
\overline{v_{k,p}(m)} = v_{k,p+1; m}.
\]
The linear independence of $v_k$, $v_{k,p}$ over 
$H=\Bbbk [\partial ]$ implies linear independence 
of the reduced words $v_{k,p;m}$.
\end{proof}

\subsection{The Anick complex for $\mathcal A_+(U(3))$}

Throughout the rest of the paper 
$\Lambda $ stands for 
the augmented algebra 
$\Lambda = \mathcal A_+(U(3))\oplus \Bbbk $
with $\varepsilon(v(n))=0$ for all $n\ge 0$.

\begin{corollary}\label{cor:U(3)AnickChains}
The Anick chains $\Lambda^{(n-1)}$ 
are of the following form:
\[
 [v(m_1)|v(m_2) | \dots |v(m_{n-1})|v(m_n)], 
\]
where $m_1,\dots, m_{n-2}\ge 2$ and either $m_{n-1}\ge 2$
or $(m_{n-1}, m_n)=(1,0)$.
\end{corollary}

We will write $[m_1|\dots |m_n]$ instead of $[v(m_1)|\dots |v(m_n)]$ for the sake of simplicity.

Let 
$\mathsf A_\bullet = 
\mathsf A_\bullet (\Lambda,\Lambda )\otimes _{\Lambda^e} \Bbbk $.
In order to calculate the differential in $\mathsf A_\bullet$
it is enough to construct a Morse matching 
in the graph $\Gamma _{\mathsf{B}_\bullet}$.
In Figures~\ref{fig:Pic-U1}--\ref{fig:Pic-U4} below,
we add extra edges corresponding to the rewriting 
of a non-reduced monomial into reduced form.
This makes it easier to track paths in the graph
$\Gamma _{\mathsf{B}_\bullet}$.

\begin{theorem}\label{thm:AnickDiff}
For $n\ge 1$, the differential  $\dcobound_{n+1}:\mathsf A_{n+1}\to \mathsf A_{n}$ is given by
\begin{multline}\label{eq:DerRegular}
    \dcobound_{n+1}[i_1|i_2|\dots |i_n|i_{n+1}]= \sum\limits^n_{j=1}
(-1)^{j}\frac{i_j(i_j-1)}{i_j+i_{j+1}-1}
      [{i_1}|{i_2}|\ldots |{i_j+i_{j+1}-1}|\ldots |{i_{n+1}}]
      \\
    +\sum\limits^n_{j=2}\sum\limits^{j-1}_{t=1}
(-1)^{j}i_t\frac{(i_j-1)(i_{j+1}-1)}{i_j+i_{j+1}-1}
 [{i_1}|\ldots |{i_t-1}|\ldots |{i_j+i_{j+1}}|\ldots | {i_{n+1}}]\\
    +\sum\limits^n_{j=2}\sum\limits^{j-1}_{t=1}
(-1)^{j}\frac{i_ji_{j+1}}{i_j+i_{j+1}-1}(i_t-1)
   [{i_1}|{i_2}|\ldots |{i_j+i_{j+1}-1}|\ldots |{i_{n+1}}],
    \end{multline}
for $i_n>1$, and 
\begin{multline}\label{eq:Der10}
\dcobound_{n+1}
[{i_1}|{i_2}|\dots |{i_{n-1}}|1|0]
= \sum\limits^{n-2}_{j=1}
(-1)^{j}\frac{i_j(i_j-1)}{i_j+i_{j+1}-1}
 [{i_1}|{i_2}|\ldots |{i_j+i_{j+1}-1}|\ldots |{i_{n-1}}| 1|0]
 \\
+ \sum\limits^{n-2}_{j=2}\sum\limits^{j-1}_{t=1}
 (-1)^{j}i_t\frac{(i_j-1)(i_{j+1}-1)}{i_j+i_{j+1}-1}
   [{i_1}|\ldots |{i_t-1}|\ldots |{i_j+i_{j+1}}|\ldots |{i_{n-1}}| 1|0]\\
+\sum\limits^{n-2}_{j=2}\sum\limits^{j-1}_{t=1}
 (-1)^{j}\frac{i_ji_{j+1}}{i_j+i_{j+1}-1}(i_t-1)
 [{i_1}|{i_2}|\ldots |{i_j+i_{j+1}-1}|\ldots |{i_{n-1}}| 1|0]\\
  +\sum\limits^{n-1}_{j=1}
  (-1)^{n}i_j
  [{i_1}|\ldots |{i_j-1}|\ldots |{i_{n-1}}|1]
+ \sum\limits^{n-1}_{j=1}
 (-1)^{n-1}(i_j-1)[{i_1}|\ldots |{i_{n-1}}|0]\\
+(-1)^n[{i_1}|\ldots |{i_{n-1}}|0].
    \end{multline}
\end{theorem}

\begin{proof}
Let us draw the segment of $\Gamma_{\mathsf B_\bullet}$ with the matched edges. First, 
draw the edges of $\Gamma_{\mathsf B_\bullet}$ 
starting at $\mathbf {b} = [i_1|\dots | i_n| i_{n+1}]$, $i_n\ge 2$,
by means of Theorem~\ref{thm:GSB-U(3)} (Fig.~\ref{fig:Pic-U1}).

Next, choose the matching for those vertices that are not 
Anick chains, i.e., those that contain 
$v(1)v(p)$ (Fig.~\ref{fig:Pic-U2}) or $v(0)v(p+1)$ (Fig.~\ref{fig:Pic-U3}) 
at $t$th position, for $t=2,\ldots, j$. 
In this way, we continue obtaining vertices of the 
same form 
until $t=1$. 

Finally, collect all Anick chains $\mathbf{d}$ obtained with the corresponding 
multiples $\Gamma_{\mathsf{A}_\bullet} (\mathbf{b},\mathbf{d})$
to get the expression for $\dcobound_{n+1}(\mathbf b)$.

For the chains ending at $[\dots |1|0]$, the proof is 
completely similar (see Fig.~\ref{fig:Pic-U4}).
\end{proof}

\begin{figure}
    \centering
      \begin{tikzpicture}[commutative diagrams/every diagram]
  \node (A) at (0,0) { $[v(i_1)|\cdots | v(i_{n+1})]$ };
  \node (B0) at (-4,-1) {$\cdots$};
  \node (B1) at (0,-2) { $[v(i_1)| \cdots |  v(i_j) v(i_{j+1})| \cdots | v(i_{n+1})  ]$ };
  \node (B2) at (4,-1) {$\cdots$\rlap{$(j=1,2,\ldots, n)$}};
  \node (C1) at (-4,-5) {$[ v(i_1)| \cdots | v(1)v(i_j+i_{j+1} -1)| \cdots v(i_{n+1}) ]$};
  \node (C2) at (4,-5) {$[v(i_1)|\cdots | v(i_j + i_{j+1}-1)| \cdots v(i_{n+1})]$};
  \node(D) at (0,-6) {$ [ v(i_1)| \cdots | v(0)v(i_j+i_{j+1})| \cdots| v(i_{n+1})]$};
  
  \path[commutative diagrams/.cd,every arrow, every label]
  (A) edge node {} (B0)
  (A) edge node {$(-1)^j$} (B1)
  (A) edge node {} (B2)
  (B1) edge node [swap]{$\frac{i_ji_{j+1}}{i_j + i_{j+1}-1}$} (C1)
  (B1) edge node [swap] {$-\frac{(i_j-1)(i_{j+1}-1)}{i_j+i_{j+1}-1}$} (D)
  (B1) edge node {$\frac{i_j(i_j-1)}{i_j + i_{j+1}-1}$} (C2)
  ;
\end{tikzpicture}
\caption{Term rewriting in the bar-complex graph for $U(3)$: generic case}
    \label{fig:Pic-U1}
\end{figure}
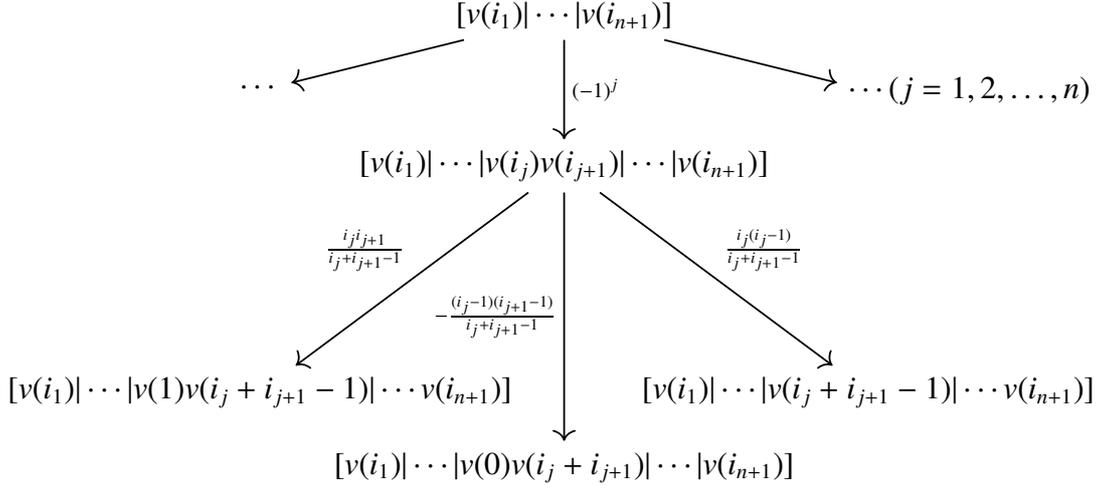

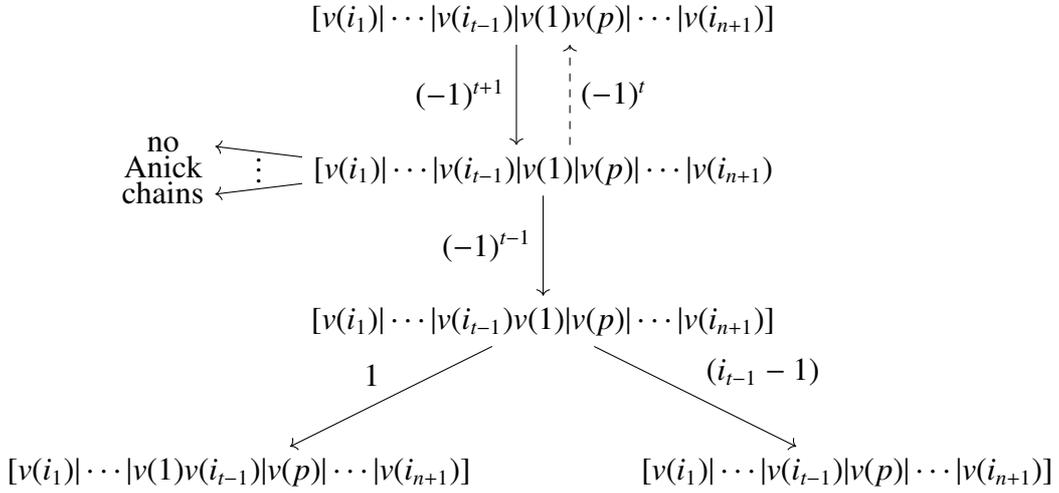
\begin{figure}
\begin{tikzpicture}[commutative diagrams/every diagram]
\node (A) at (0,0) {$[v(i_1)|\cdots | v(i_{t-1})|v(1)v(p)|\cdots| v(i_{n+1})]$};
\node (B) at (0,-2) {$[v(i_1)|\cdots| v(i_{t-1})|v(1)|v(p)| \cdots| v(i_{n+1})$};
\node(B0) at(-5,-2) {$\substack{ \mbox{no} \\ \mbox{Anick} \\ \mbox{chains} }$};
\node (C) at (0,-4) {$[v(i_1)|\cdots| v(i_{t-1})v(1)|v(p)|\cdots |v(i_{n+1})]$};
\node (D1) at (-4, -6) {$[v(i_1)|\cdots | v(1)v(i_{t-1})| v(p)|\cdots| v(i_{n+1})]$};
\node (D2) at (4,-6) {$[v(i_1)|\cdots | v(i_{t-1})| v(p)| \cdots | v(i_{n+1})]$};

\path[->]
  ([xshift=-10pt]A.south) edge node[left] {$(-1)^{t+1}$} ([xshift=-10pt]B.north);
\path[->,dashed]
 ([xshift=10pt]B.north) edge node[right]{$(-1)^{t}$} ([xshift=10pt]A.south);
\path[->]
 ([yshift=5pt]B.west) edge node[yshift=-4pt] {$\vdots$} ([yshift=9pt]B0.east);
\path[->]
 ([yshift=-5pt]B.west) edge node {} ([yshift=-9pt]B0.east);
\path[->]
 (B) edge node [left] {$(-1)^{t-1}$} (C);
\path[->]
  (C) edge node [above left] {$1$} (D1)
  (C) edge node [above right] {$(i_{t-1}-1)$} (D2);
\end{tikzpicture}
    \caption{A segment of the Morse matching graph for $U(3)$}
    \label{fig:Pic-U2}
\end{figure}

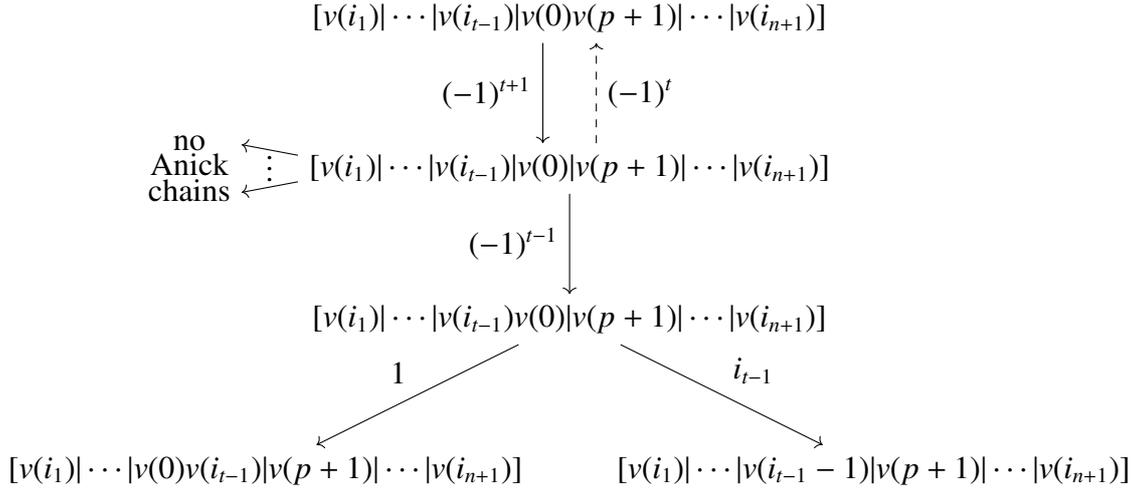
\begin{figure}
   \begin{tikzpicture}[commutative diagrams/every diagram]
\node (A) at (0,0) {$[v(i_1)|\cdots | v(i_{t-1})|v(0)v(p+1)|\cdots| v(i_{n+1})]$};
\node (B) at (0,-2) {$[v(i_1)|\cdots| v(i_{t-1})|v(0)|v(p+1)| \cdots| v(i_{n+1})]$};
\node(B0) at(-5,-2) {$\substack{ \mbox{no} \\ \mbox{Anick} \\ \mbox{chains} }$};
\node (C) at (0,-4) {$[v(i_1)|\cdots| v(i_{t-1})v(0)|v(p+1)|\cdots |v(i_{n+1})]$};
\node (D1) at (-4, -6) {$[v(i_1)|\cdots | v(0)v(i_{t-1})| v(p+1)|\cdots| v(i_{n+1})]$};
\node (D2) at (4,-6) {$[v(i_1)|\cdots | v(i_{t-1}-1)| v(p+1)| \cdots | v(i_{n+1})]$};

\path[->]
  ([xshift=-10pt]A.south) edge node[left] {$(-1)^{t+1}$} ([xshift=-10pt]B.north);
\path[->,dashed]
 ([xshift=10pt]B.north) edge node[right]{$(-1)^{t}$} ([xshift=10pt]A.south);
\path[->]
 ([yshift=5pt]B.west) edge node[yshift=-4pt] {$\vdots$} ([yshift=9pt]B0.east);
\path[->]
 ([yshift=-5pt]B.west) edge node {} ([yshift=-9pt]B0.east);
\path[->]
 (B) edge node [left] {$(-1)^{t-1}$} (C);
\path[->]
  (C) edge node [above left] {$1$} (D1)
  (C) edge node [above right] {${i_{t-1}}$} (D2);
\end{tikzpicture}
    \caption{A segment of the Morse matching graph for $U(3)$}
    \label{fig:Pic-U3}
\end{figure}

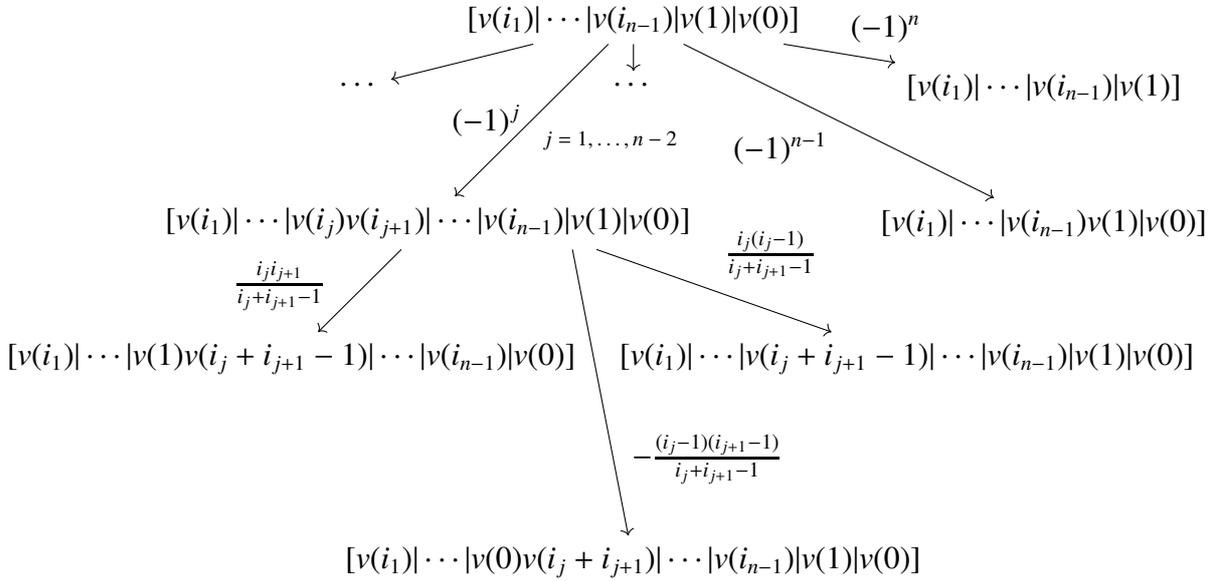
\begin{figure}
 \begin{tikzpicture}[scale=0.9,commutative diagrams/every diagram]
\node (A) at (0,0) {$[v(i_1)|\cdots|v(i_{n-1})|v(1)|v(0)]$};

\node (B1) at (-4,-1) {$\cdots$};
\node (B2) at (0,-1) {$\cdots$};
\node (B3) at (6,-1) {$[v(i_1)|\cdots| v(i_{n-1})|v(1)]$};

\node (C1) at (-3,-3) {$[v(i_1)|\cdots | v(i_j) v(i_{j+1})|\cdots | v(i_{n-1})|v(1)|v(0)]$};
\node (C2) at (6,-3) {$[v(i_1)|\cdots | v(i_{n-1})v(1)|v(0)]$};

\node (D1) at (-5,-5) {$[v(i_1)|\cdots | v(1)v(i_j + i_{j+1}-1)|\cdots| v(i_{n-1})|v(0)]$};
\node (D2) at (4,-5) {$[v(i_1)|\cdots | v(i_j+i_{j+1}-1)|\cdots|v(i_{n-1})|v(1)|v(0)]$};

\node (E) at (0,-8) {$[v(i_1)|\cdots| v(0)v(i_j+i_{j+1})|\cdots|v(i_{n-1})|v(1)|v(0)]$};

\path[->]
 (A) edge node {} (B1)
 (A) edge node[left] {$(-1)^j$} node[below right] {\small{\tiny{$j=1,\ldots, n-2$}}} (C1)
 (A) edge node {} (B2)
 (A) edge node [below left] {$(-1)^{n-1}$} (C2)
 (A) edge node [above right] {$(-1)^n$} (B3)
 ;

\path[->]
 (C1) edge node [xshift = -30pt, yshift = 2pt] {$\frac{i_j i_{j+1}}{i_j + i_{j+1}-1}$} (D1)
 ([xshift = 70pt]C1.south) edge node [above right] {$\frac{i_j(i_j-1)}{i_j+i_{j+1}-1}$} (D2)
 ([xshift = 60pt]C1.south) edge node [yshift=-25pt, xshift = 40pt] {$-\frac{(i_j-1)(i_{j+1}-1)}{i_j + i_{j+1}-1}$} (E);
\end{tikzpicture}
    \caption{Term rewriting in the bar-complex graph for $U(3)$: chains ending with $[\dots |1|0]$}
    \label{fig:Pic-U4}
\end{figure}

\begin{example}
Let $[n|m]$ be an Anick 1-chain in $\Lambda $
which is not equal to $[1|0]$. 
Then 
\[
\dcobound_2([n|m]) = -\dfrac{n(n-1)}{n+m-1} [n+m-1].
\]
For the remaining 1-chain we have 
\[
\dcobound_2([1|0]) = -[0].
\]
\end{example}

\begin{example}
Let $[n|m|p]$ be an Anick 2-chain in $\Lambda $
with $m\ge 2$. 
Then $\dcobound_3([n|m|p])$ is equal to
\[
\dfrac{m(p(n-1)+m-1)}{m+p-1} [n|m+p-1]
-\dfrac{n(n-1)}{n+m-1} [n+m-1|p] - \dfrac{n(m-1)(p-1)}{m+p-1} [n-1|m+p].
\]
(For $n=2$, the last summand is missing.)
In a similar way, 
\[
\dcobound_3([n|1|0]) = (2-n)[n|0]+n[n-1|1], \quad n\ge 2.
\]
\end{example}

Recall that $\Lambda $ carries a derivation $\partial $
given by $\partial (v(n)) = nv(n-1)$.
Let us describe the corresponding map $\tilde \partial_\bullet$
on $\mathsf A_\bullet $. 
According to the general scheme \eqref{eq:DerivTilda}, 
in order to calculate 
$\tilde \partial_n ([m_1|\dots |m_n])$
one has to evaluate 
$\mathrm g_n([m_1|\dots |m_n]) \in \mathsf B_n$, 
apply $\partial_n$ by \eqref{eq:Partial_on_B}, 
and remove all those summands that are not Anick chains
(i.e., apply $\mathrm f_n$).

Note that if $[m_1|\dots |m_{n}]\in \mathsf B_{n}$
is not an Anick chain then the expression for 
$\partial_n ([m_1|\dots |m_n])$ does not contain 
Anick chains except for the case when 
$m_{n-1}=m_n=1$, $m_1,\dots, m_{n-2}\ge 2$.
It can be seen from the proof of Theorem \ref{thm:AnickDiff}
that such a summand
$\mathbf u = [m_1|\dots |m_{n-2}|1|1] \in \mathsf B_n$
may appear only in 
$\mathrm g_n([m_1|\dots |m_{n-2}|2|0])$, 
but even in this case the coefficient at 
$\mathbf u$ 
is equal to zero
since $v(2)v(0)= v(0)v(2) + 2v(1)$, does not 
contain $v(1)v(1)$.

Hence, the calculation of $\tilde\partial _n$ on $\mathsf A_n$
becomes quite simple: one has to differentiate 
an Anick chain $[m_1|\dots |m_{n}]$ as it was an 
element of $\mathsf B_n$, 
and then remove all those 
summands that are not Anick chains. 
For example, 
$\tilde\partial _2([2|p]) = p[2|p-1]$,
$\tilde\partial _3([n|2|1]) = n[n-1|2|1] + [n|2|0]$ if $n>2$,
$\tilde\partial _n([2|2|\dots |2|1|0]) = 0$,
etc.

For every $n\ge 1$,
the space $\mathsf A_n$ has the following  grading:
\[
\mathsf A_n^{(d)} = \{\mathbf a \in \mathsf A_n \mid \deg \mathbf a = d \},
\]
where $\deg [i_1|\ldots |i_n] = i_1+\cdots + i_n$.
Theorem~\ref{thm:AnickDiff} shows $\dcobound_n$ to be 
a degree $-1$ map: 
$\dcobound_n(\mathsf A_n) \subseteq \mathsf A_{n-1}^{(d-1)}$.

The space 
$K_n=\Ker\tilde\partial_n\subseteq \mathsf A_n$ 
inherits the  grading from $\mathsf A_n$:
$K_n^{(d)} = K_n\cap \mathsf A_n^{(d)}$.

\subsection{Low-dimensional Hochschild conformal cohomology
of $U(3)$}\label{subsec:1-3dim}
Let us explicitly calculate the dimensions of the
first, second, and third 
Hochschild cohomology groups for the 
associative conformal algebra $U(3)$ with coefficients in the scalar 
module~$\Bbbk $. 

Denote $K_n = \Ker \tilde \partial_n \subseteq \mathsf A_n $
for $n\ge 1$. According to the general principle
from Section~\ref{subsec:ReducedComplex}, 
$\Homol ^n (U(3),\Bbbk )\simeq (\Ker \tilde \dcobound_n/\im \tilde \dcobound_{n+1} )^*$, 
where $\tilde\dcobound_n $ 
is the restriction of 
$\dcobound_n: \mathsf A_n \to \mathsf A_{n-1}$
onto $K_n$.

For example, $\dim K_1=1$, the basis is $[0]$,
and $\dcobound_2([1|0]) = -[0]$, 
where $[1|0]\in K_2$. Hence, 
$\Homol^1(U(3),\Bbbk ) = 0$.

An arbitrary element $\mathbf a \in \mathsf A_n $
may be uniquely written as 
\[
\mathbf a = [\mathbf a_0|0]+[\mathbf a_1|1]+\dots,
\]
where $\mathbf a_i \in \mathsf A_{n-1}$.
%
% a_i : all chans like [*** | i]
%

\begin{lemma}\label{lem:Zero-Principle}
Every element $\mathbf a \in K_n$ 
is completely determined by $\mathbf a_0\in \mathsf A_{n-1}$
which satisfies the following properties:
{\em (1)} $(\mathbf a_0)_0 = 0$; 
{\em (2)} $(\tilde \partial _{n-1} \mathbf a_0)_1 = 0$.
\end{lemma}

\begin{proof}
Suppose $\mathbf b = \tilde\partial_n \mathbf a$, 
$\mathbf b = [\mathbf b_0 | 0] + [\mathbf b_1 | 1] +\dots $.
Then
\[
[\mathbf b_i | i] = \mathrm f_n \big ([\tilde \partial_{n-1} (\mathbf a_i)|i] 
 + (i+1)[\mathbf a_{i+1} | i]\big ) , 
 \quad i=0,1,\dots .
\]
Note that $\mathrm f_n ([\mathbf a_{i+1}|i]) 
= [\mathbf a_{i+1}|i]$ for $i\ge 0$.

Then $\mathbf a \in K_n$ if and only if $\mathbf b_i=0$ 
for all $i\ge 0$. Hence, $\mathbf a_{i+1}$ is completely
determined by 
$\mathrm f_n ([\tilde \partial_{n-1} (\mathbf a_i)|i] )$.
The only problem emerges for $i=0$: 
 $\mathrm f_n ([\tilde \partial_{n-1} (\mathbf a_0)|0])$
 may contain a chain of the form $[\cdots |1|0]$, 
 but $\mathbf a_1$ may not contain a chain ending with 
 $[\cdots |1]$. Therefore, $\mathbf a$ 
 is completely defined by $\mathbf a_0$ such that 
 $[\mathbf a_0 |0] $ is a combination of chains
 (i.e., $\mathbf a_0 \in \mathsf A_{n-1}$
 does not contain zeros) and 
 $\tilde \partial_{n-1} (\mathbf a_0)$
 does not contain units.
\end{proof}

Let $\mathsf A_n^{\circ }$ stand for the space 
spanned by {\em regular} chains, i.e., 
those that have all components $\ge 2$.
Similarly, denote 
$\mathsf A_n^{(d)\circ } = \mathsf A_n^{\circ }\cap \mathsf A_n^{(d)}$.

\begin{corollary}\label{cor:Kn-structure}
For $n\ge 3$ we have 
$K_n \simeq  \mathsf A_{n-2}^{(d-1)\circ }\oplus 
\mathsf A_{n-2}^{(d-3)\circ } \oplus \mathsf A_{n-2}^{(d-4)\circ } \oplus \dots $.
\end{corollary}

\begin{proof}
Assume $n\ge 3$ and $\mathbf v = \mathbf a_0$ 
for $\mathbf a \in K_n$ meets the conditions (1), (2) 
of Lemma~\ref{lem:Zero-Principle}.
Then $\mathbf v\in \mathsf A_{n-1}$
may be presented in the same form 
\[
\mathbf v = [\mathbf v_1|1]+[\mathbf v_2|2]+\dots,
\]
where $\mathbf v_i\in \mathsf A_{n-2}^{(d-i)\circ }$.

The condition $(\tilde\partial_{n-1}\mathbf v)_1=0$
is equivalent to
\[
\mathrm f_{n-1}\big ([\tilde\partial_{n-2}(\mathbf v_1)|1]\big )
+2[\mathbf v_2|1] = 0.
\]
Hence, $\mathbf v_1$ uniquely 
determines $\mathbf v_2$ and,
together with $\mathbf v_k$, $k\ge 3$, 
uniquely define $\mathbf a \in K_n$.
\end{proof}

\begin{corollary}
The linear basis of $K_2$ consists of 
\[
\mathbf e_1=[1|0], \ \mathbf e_3=[3|0]-3[2|1], \ \dots, \ 
\mathbf e_d =\sum\limits_{s=0}^{d-2} (-1)^s \binom{d}{s} [d-s|s], 
\]
for $d\ge 3$.
\end{corollary}

\begin{proposition}[c.f. \cite{AlKol-JMP}]\label{prop:H2-U(3)}
$\dim \Homol^2(U(3),\Bbbk )=1$.
\end{proposition}

\begin{proof}
Note that $\mathbf e_d$ 
for $d\ge 3$ belong to $\Ker \dcobound_2$.
One may either apply Theorem~\ref{thm:AnickDiff} or just note that 
 $\dcobound_2(K^{(d)}_2)\subseteq K_1^{(d-1)}$ 
and $K_1^{(d-1)}=0$ for $d\ne 1$.

On the other hand, for every $d\ge 3$
%otherwise, \partial has [...1]
there exists $\mathbf f_{d+1} \in K_3^{(d+2)}$ such that 
\[
\mathbf f_{d+1} = [2|d|0] + [\cdots |1] +\dots.
\]
Then $\dcobound_3(\mathbf f_{d+1}) 
 = -\frac{2}{d+1}[d+1|0] + [\cdots |1]+\dots \in K_2^{(d+1)}$, 
so $\mathbf e_{d+1} \in \im \tilde \dcobound_3$.
Therefore, it remains to compare $\Ker \tilde \dcobound_2$ and 
$\im \tilde \dcobound_3$
in $K_2^{(3)}$. The kernel is 1-dimensional, but the image is zero. 
Indeed, the only regular chain in $\mathsf A_2^{(4)}$ is $[2|2]$, so 
$\dim K_3^{(4)}=1$ and the basis (e.g., recovered by Lemma~\ref{lem:Zero-Principle}) is
\[
\mathbf f_3 = [2|2|0] - \dfrac{2}{3} [3|1|0].
\]
It is straightforward to compute that 
$\dcobound_3(\mathbf f_3)=0$, so $\im \tilde\dcobound_3 =0$.
\end{proof}

\begin{proposition}
$\dim \Homol^3(U(3),\Bbbk )=1$.
\end{proposition}

\begin{proof}
Let $\tilde\dcobound_n^{d}$ stand for the 
differential $K_n^{(d)} \to K_{n-1}^{(d-1)}$. 
By Corollary~\ref{cor:Kn-structure} $\dim K_3^{(d)}=d-3$
for $d\ge 4$. On the other hand, 
$\im \tilde \dcobound_3^{(d)} = K_2^{(d-1)}$ for $d-1\ge 4$, as we have seen in the proof 
of Proposition~\ref{prop:H2-U(3)}. 
Hence, $\dim\Ker \tilde\dcobound_3^{(d)} = d-4$ 
for 
$d\ge 5$.

Suppose $d\ge 5$. Choose the following elements 
$\mathbf a^{(i)}\in K_4^{(d+1)}$, $i=2,\ldots, d-4$, 
by defining their zero components 
\[
\mathbf a_0^{(i)} = [2|i|d-1-i].
\]
One more element $\mathbf a^{(d-3)} \in K_4^{(d+1)}$
is defined by 
\[
\mathbf a^{(d-3)}_0 = [2|d-3|2] - \dfrac{2}{d-2} [2|d-2|1].
\]
Let us calculate the leading (in the lexicographic sense) terms of $\tilde \dcobound_4^{(d+1)} (\mathrm a^{(i)})$
by Theorem~\ref{thm:AnickDiff}:
\[
\tilde \dcobound_4^{(d+1)} (\mathbf a^{(i)})
 = \alpha_i [i+1|d-i-1|0] +\dots ,\quad i=2,3,\dots, d-3,
\]
where $\alpha_i\in \Bbbk $, $\alpha _i\ne 0$.
Hence, $\dim\im \tilde\dcobound_4^{(d+1)} = d-4 = \dim\Ker \tilde \dcobound_3^{(d)}$.

Therefore, nontrivial cohomology may appear in $K_3^{(d)}$
for $d<5$ only. Recall that $K_3^{(4)}$ is spanned 
by $2[3|1|0]-3[2|2|0]$ which is proportional 
to the image of $[2|2|1|0]\in K_4^{(5)}$ under $\tilde \dcobound_4$.
For $d=3$, we have only $[2|1|0]\in \Ker\tilde \dcobound_3^{(3)}$ which may not be an image 
of $\tilde\dcobound_4^{(4)}$ since there are no 
Anick chains of degree 4 in $\Lambda ^{(3)}$.
\end{proof}

\subsection{Higher Hochschild cofomology of $U(3)$}\label{4.3}
% New version
There is a filtration 
\[
\mathsf A_{n,0}\supset \mathsf A_{n,1} \supset \mathsf A_{n,2} \supset \dots,
\]
where 
\[
\mathsf A_{n,k} = \mathrm{span}\big \{[i_1|\ldots |i_n] \in \Lambda^{(n-1)} 
 \mid i_n\ge k \big \}.
\]
It follows from Theorem~\ref{thm:AnickDiff} that 
$\dcobound_n (\mathsf A_{n,k}) \subseteq \mathsf A_{n-1,k}$.

\begin{proposition}\label{prop:NonRestricted-Regular}
Let $n\ge 2$, 
$\mathbf u\in \mathsf A_{n,2}$, 
and $\dcobound_n\mathbf u =0$. 
Then there exists $\mathbf w \in \mathsf A_{n+1,2}$
such that $\dcobound_{n+1}\mathbf w = \mathbf u$.
\end{proposition}

\begin{proof}
We may suppose that $\mathbf u \in \mathsf A_n^{(d)}$, $d\ge 5$
(since $[2|2]$ is the only regular chain of degree 4).
For $n=2$, the proof is straightforward. 
Indeed, assume the converse and choose
\begin{equation}\label{eq:Ker-not-Cobound}
\mathbf u = \sum\limits_{i=0}^{d-4} \alpha_i [d-i-2|i+2] 
\in 
\big (\Ker\dcobound_2\cap \mathsf A_{2,2}^{(d)} \big)
\setminus \dcobound_3(\mathsf A_{3,2})
\end{equation}
with minimal $l$ such that $\alpha _l\ne 0$.
Note that $l<d-4$ since 
$\dcobound_2([2|d-2])\ne 0$.
Then 
we may build
\[
\mathbf u' = \mathbf u - \gamma \dcobound_{3}[d-l-1|2|i+2]
\]
for $\gamma = -((d-l-1)(d-l-2)/(d-l))^{-1}$
which has the same property as $\mathbf u$
in \eqref{eq:Ker-not-Cobound}, 
but its presentation has nonzero coefficients for $i>l$, 
which is a contradiction.

Suppose $n>2$ and assume the statement is true for 
$\mathbf v \in \mathsf A_{m,2}$ with $m<n$.
Then 
\[
\mathbf u = [\mathbf u_2|2]+\mathbf v, \quad \mathbf v\in \mathsf A_{n,3}, 
\]
and $\dcobound_n\mathbf u=0$ implies $\dcobound_{n-1}\mathbf u_2=0$
since $\dcobound_n\mathbf v \in \mathsf A_{n-1,3}$.
Hence there exists 
$\mathbf w_2 \in \mathsf A_{n,2}$
such that $\dcobound_n \mathbf w_2 = \mathbf u_2$. 
Consider
$\mathbf u' = \mathbf u - \dcobound_{n+1}[\mathbf w_2|2]$.
Since 
$\dcobound_{n+1}[\mathbf w_2|2]\in [\dcobound_n\mathbf w_2|2]
+\mathsf A_{n,3}$
by Theorem~\ref{thm:AnickDiff}, we obtain 
\[
\mathbf u' \in  \Ker\dcobound_n\cap \mathbf A_{n,3}.
\]

In a similar way, we may find $\mathbf w_3 \in \mathsf A_{n,2}$
such that 
$\mathbf u' - \dcobound_{n+1}[\mathbf w_3|3] \in \mathsf A_{n,4}$, 
and so on. Since $\mathsf A_{n,k}^{(d)}=0$ for sufficiently 
large $k$, the result follows.
\end{proof}

\begin{corollary}\label{cor:NonRestr-1}
Let $\mathbf u\in \mathsf A_{n,1}$, $n\ge 3$, and $\dcobound_n \mathbf u=0$.
Then there exists $\mathbf w \in \mathsf A_{n+1,1}$ such that 
$\dcobound_{n+1}\mathbf w  = \mathbf u$.
\end{corollary}

\begin{proof}
Let $\mathbf u = [\mathbf u_1|1] + \mathbf u'$, $\mathbf u' \in \mathsf A_{n,2}$.
Then $\dcobound_n\mathbf u=0$ implies $\dcobound_{n-1}\mathbf u_1=0$.
Since $\mathbf u_1\in \mathsf A_{n-1,2}$ and
$n-1\ge 2$, there exists $\mathbf w_1\in \mathsf A_{n,2}$
such that $\dcobound_n\mathbf w_1 = \mathbf u_1$.
Then 
$\mathbf u - \dcobound_{n+1}[\mathbf w_1|1] \in \mathsf A_{n,2}$
may be presented as $\dcobound_{n+1}\mathbf w$
for an appropriate $\mathbf w \in \mathsf A_{n+1,2}$.
\end{proof}

\begin{corollary}\label{cor:NonRestr-0}
Let $\mathbf u\in \mathsf A_{n}$, $n\ge 4$, 
and $\dcobound_n \mathbf u=0$.
Then there exists $\mathbf w \in \mathsf A_{n+1}$ 
such that 
$\dcobound_{n+1}\mathbf w  = \mathbf u$.
\end{corollary}

\begin{proof}
If $\mathbf u\in \mathsf A_{n,1}$ then we are done. 
Hence, consider 
$\mathbf u = [\mathbf u_0|0] +\mathbf u'$, $\mathbf u'\in \mathsf A_{n,1}$.
Then $\mathbf u_0\in \mathsf A_{n-1,1}$ and $\dcobound_{n-1}\mathbf u_0 = 0$.
By Corollary \ref{cor:NonRestr-1} there exists $\mathbf w_0\in \mathsf A_{n,1}$
such that $\dcobound_n\mathbf w_0 = \mathbf u_0$.
Then 
$\mathbf u - \dcobound_{n+1}[\mathbf w_0|0] \in \mathsf A_{n,1}$
and thus $\mathbf u = \dcobound_{n+1}\mathbf w$
for an appropriate $\mathbf w \in \mathsf A_{n+1}$.
\end{proof}

\begin{theorem}
$\Homol^n(U(3),\Bbbk )=0$ for $n\ge 4$.
\end{theorem}

\begin{proof}
Corollary~\ref{cor:NonRestr-0} implies 
the cohomology groups of the non-restricted 
complex $\widetilde\C^\bullet$  of the conformal algebra 
$U(3)$ with coefficients in $M=\Bbbk $ are trivial 
for $n\ge 4$.
The short exact sequence of complexes 
\[
0\to \partial\widetilde \C^\bullet \to \widetilde\C^\bullet \to 
\C^\bullet \to 0
\]
leads to the long exact sequence of cohomology groups
\[
\begin{aligned}
\dots & \to \Homol^n(\partial\widetilde \C^\bullet)
        \to \Homol^n(\widetilde \C^\bullet)
        \to \Homol (\C^\bullet) \\
& \to \Homol^{n+1}(\partial\widetilde \C^\bullet) \to \dots .
\end{aligned}
\]
By \cite[Proposition~2.1]{BKV}
$\partial\widetilde \C^\bullet\simeq \widetilde \C^\bullet$
in degrees $n\ge 1$, so 
$\Homol^n(\widetilde \C^\bullet) = \Homol^n(\partial\widetilde \C^\bullet) = 0$
for $n\ge 4$.
Hence, the restricted complex $\C^\bullet$
has zero cohomologies for $n\ge 4$.
\end{proof}

\subsection{Final remarks}
We have shown that the Hochschild cohomology groups 
for the universal associative conformal envelope 
$U(3) = U(\Vir; N=3)$ with coefficients in the scalar module 
are the same as for the Virasoro Lie conformal algebra 
$\Vir $. However, this is not true for the cohomologies 
with coefficients in a non-trivial irreducible 
module. 

Recall that all irreducible modules over $\Vir$ 
are of the form $M(\alpha,\Delta)$, $\Delta\ne 0$ \cite{ChengKac}, 
see also Example~\ref{exmp:Vir-modules}.
It was found in \cite{BKV} that 
\[
\dim \Homol^1(\Vir, M(0,\Delta))
 = \begin{cases}
 2, & \Delta  = 1, \\
 1, & \Delta=-1,0, \\
 0, & \text{otherwise}.
 \end{cases}
\]

All these representations may be extended to the 
representations of $U(3)$, so, in particular, $M(0,\Delta)$
is a left $U(3)$-module for every scalar $\Delta $. 
(Note that the envelope $U(2)$ inherits only those representations 
that correspond to $\Delta =0,1$.)

Let us compute 
$\Homol^1(U(3), M)$ for $M=M(0,\Delta)$, $\Delta \in \Bbbk $,
in a straightforward way. 

Let $\varphi \in \Z^1(U(3),M)$. This is a conformal module map, 
i.e., a $\Bbbk [\partial ]$-linear map from $U(3)$ to $M$
such that 
\[
\varphi (a\oo{\lambda } b) = a\oo{\lambda } \varphi(b),
\quad a,b\in U(3).
\]
The value $\varphi(v)=f(\partial )u$, $f\in \Bbbk [\partial ]$,
on the generator of $U(3)$
determines the cocycle completely. 
In order to preserve the defining relation \eqref{eq:vir-comm}, 
the polynomial $f(\partial )$ has to meet the following conditions:
\[
\begin{gathered}
2v\oo{1} f(\partial )u -\partial (v\oo{2} f(\partial)u) = 2f(\partial )u, \\
v\oo{n} f(\partial ) u = 0, \quad n\ge 3.
\end{gathered}
\]
The sesqui-linearity implies 
$v\oo{n} f(\partial )u = (\partial f^{(n)}(\partial ) + n\Delta f^{(n-1)}(\partial ))u$, for $n\ge 0$.
Hence, for $\Delta \ne 0$
we have $f^{(2)}=f^{(3)}= \dots =0$ from the locality condition, 
so $f(\partial ) = c_0 +\partial c_1$.
The latter meets 
\[
(\partial f' +\Delta f) - \partial \Delta f' = f,
\]
which turns into $(\Delta-1)c_0 = 0$. 
Hence, $\dim \Z^1(U(3),M) = 1$ for $\Delta \ne 1,0$
and $\dim \Z^1(U(3),M) = 2$ for $\Delta =1$.

If $\Delta=0$ then we obtain $f^{(3)}=f^{(4)}=\dots =0$
from the locality condition, thus $f = c_0+c_1\partial +c_2\partial^2$
in this case. Then the defining relation turns into
\[
2\partial f' -\partial^2 f'' = 2f,
\]
which implies the only condition $c_0=0$.
Therefore, $\dim \Z^1(U(3),M)=2$ for $\Delta =0$.

The space of coboundaries $\B^1(U(3),M)$ consists of those $\varphi \in \Z^1(U(3),M)$ that can be presented as 
$\varphi (a) = \{a\oo{0} u\}= (a\oo{-\partial } u)$, 
$a\in U(3)$. 
For $a=v$, we obtain $\varphi(v) = f(\partial )u = (-\partial\Delta +\partial )u $.
Hence, $\dim \B^1(U(3),M)=1$ for $\Delta \ne 1$ and $\B^1(U(3),M)=0$
for $\Delta= 1$.

As a result, 
for the 1st Hochschild cohomologies of $U(3)$ we have 
\[
\dim \Homol^1(U(3), M(0,\Delta))
= \begin{cases}
2, & \Delta  = 1, \\
1, & \Delta=0, \\
0, & \text{otherwise}. 
  \end{cases}
\]
Therefore, there is a difference between
cohomologies of $\Vir $ and its 
universal associative envelope $U(3)$.
An interesting task is to calculate 
all Hochschild cohomology groups
for $U(3)$ with coefficients in 
$M(\alpha,\Delta) $.

\end{document}